\documentclass{amsart}

\usepackage{amssymb}
\usepackage[all]{xy}

\newtheorem{thm}{Theorem}

\newtheorem{lem}[thm]{Lemma}
\newtheorem{cor}[thm]{Corollary}

\newtheorem{prop}[thm]{Proposition}

\newtheorem{conj}[thm]{Conjecture}
   
\theoremstyle{definition}
\newtheorem{defn}[thm]{Definition}
\newtheorem{condition}[thm]{Condition}
\newtheorem{say}[thm]{}
\newtheorem{exmp}[thm]{Example}

\newtheorem{ques}[thm]{Question}    

\newtheorem{rem}[thm]{Remark}          

\newtheorem*{ack}{Acknowledgments}      

\newtheorem{defn-thm}[thm]{Definition--Theorem}  
\newtheorem{defn-lem}[thm]{Definition--Lemma}  

\theoremstyle{remark}
\newtheorem{claim}[thm]{Claim}

\let \cedilla =\c
\renewcommand{\c}[0]{{\mathbb C}}  

\renewcommand{\o}[0]{{\mathcal O}} 
\newcommand{\z}[0]{{\mathbb Z}}
\newcommand{\n}[0]{{\mathbb N}}

\renewcommand{\a}[0]{{\mathbb A}}

\newcommand{\p}[0]{{\mathbb P}}
\newcommand{\f}[0]{{\mathbb F}}
\newcommand{\q}[0]{{\mathbb Q}}

\newcommand{\qtq}[1]{\quad\mbox{#1}\quad}
\newcommand{\spec}[0]{\operatorname{Spec}}
\newcommand{\pic}[0]{\operatorname{Pic}}
\newcommand{\pics}[0]{\operatorname{\mathbf{Pic}}}
\newcommand{\picl}[0]{\operatorname{Pic}^{\rm loc}}
\newcommand{\picel}[0]{\operatorname{Pic}^{\rm et-loc}}
\newcommand{\pical}[0]{\operatorname{Pic}^{\rm an-loc}}
\newcommand{\pica}[0]{\operatorname{Pic}^{\rm an}}
\newcommand{\picls}[0]{\operatorname{\mathbf{Pic}}^{\rm loc}}
\newcommand{\piclf}[0]{\operatorname{{\it Pic}}^{\rm loc}}
\newcommand{\piclo}[0]{\operatorname{\mathbf{Pic}}^{\rm loc-\circ}}
\newcommand{\piclt}[0]{\operatorname{\mathbf{Pic}}^{\rm loc-\tau}}

\newcommand{\pico}[0]{\operatorname{\mathbf{Pic}}^{\circ}}
\newcommand{\pict}[0]{\operatorname{\mathbf{Pic}}^{\rm \tau}}

\newcommand{\nsl}[0]{\operatorname{NS}^{\rm loc}}

\newcommand{\supp}[0]{\operatorname{Supp}}    
\newcommand{\red}[0]{\operatorname{red}}    
\newcommand{\codim}[0]{\operatorname{codim}}    
\newcommand{\im}[0]{\operatorname{im}}

\newcommand{\coker}[0]{\operatorname{coker}}    
    
\newcommand{\Hom}[0]{\operatorname{Hom}}

\newcommand{\chr}[0]{\operatorname{char}}

\newcommand{\len}[0]{\operatorname{length}}

\newcommand{\rdown}[1]{\lfloor{#1}\rfloor}

\newcommand{\onto}[0]{\twoheadrightarrow}

\newcommand{\obs}[0]{\operatorname{obs}}

\newcommand{\depth}[0]{\operatorname{depth}} 
\newcommand{\tsum}[0]{\textstyle{\sum}}

\newcommand{\shom}[0]{\operatorname{\mathcal{H}\!\it{om}}}

\newcommand{\defor}[0]{\operatorname{Def}} 
\newcommand{\deforb}[0]{\operatorname{\mathbf{Def}}}

\def\into{\DOTSB\lhook\joinrel\to}

\def\loccoh#1.#2.#3.#4.{H^{#1}_{#2}(#3,#4)}

\DeclareMathAlphabet{\mathchanc}{OT1}{pzc}%
                                {m}{it}

\newcommand{\gm}[0]{{\mathbb G}_m}

\usepackage[all]{xy}\xyoption{dvips}

\newcommand{\norm}[0]{\operatorname{norm}}

\newcommand{\link}[0]{\operatorname{link}}

\begin{document}
\bibliographystyle{amsalpha}


\title{Maps between   local Picard groups}
\author{J\'anos Koll\'ar}

\maketitle

Let $X$ be a scheme and $x\in X$ a  point. 
The {\it local Picard group}
of $X$ at $x$, denoted by
$\picl(x,X) $, is the Picard group of
the {\it punctured neighborhood}  $\spec_X\o_{x,X} \setminus \{x\}$.
Our aim is to study the pull-back map on the local Picard group
in two situations.

\begin{ques}[Normalization] \label{main.ques.1}
Let $\pi:\bar X\to X$ denote the normalization
and $\bar x_i\in \bar X$ the  preimages of $x$. What is the kernel of the
pull-back map
$$
\pi^* : \picl(x,X)\to \tsum_i \picl(\bar x_i,\bar X)?
$$
\end{ques}

\begin{ques}[Restriction to a divisor] \label{main.ques.2}
Let  $x\in D\subset X$ be an effective
 Cartier divisor.  What is the kernel of the
restriction map
$$
r^X_D : \picl(x,X)\to\picl(x,D)?
$$
\end{ques}

In both cases we are interested in conditions that guarantee
that the kernels of these maps are ``small.'' 
Here ``small'' can mean trivial, or finite or ``naturally'' a subgroup of
a linear algebraic group.

Grothendieck's local Lefschetz-type theorem
\cite[XI.3.16]{sga2} is the first major result on
Question \ref{main.ques.2}. He proves that  $r^X_D$ is injective
if $\depth_xX\geq 4$. 
  Conjecture 1.2 in  \cite{k-gl1} asserts
that this can be relaxed to $\depth_xX\geq 3$ and $\dim X\geq 4$. 
This was proved in \cite{k-gl1} when $X$ is log canonical
and in \cite{bha-dej} when $X$ is normal.

I was led to  Question \ref{main.ques.1}
 while investigating stable varieties and their
moduli. The method of \cite{kk-singbook} studies first the normalization 
$\bar X$ of a stable variety $X$ and then descends the information from
$\bar X$ to $X$. It is especially important to know when the canonical divisor
 of $X$ is $\q$-Cartier; this was settled in \cite{k-source}.
For  semi-log canonical  varieties Question \ref{main.ques.1} 
is answered in \cite[Sec.5.7]{kk-singbook}.

In many applications, the natural assumption is that 
 $X$ satisfies Serre's condition $S_2$ only. 
To understand this case, 
one needs to know that, as in the global situation, the 
local Picard group is not just a group but a  group scheme. The local
 definition is  more
involved than the projective one and the natural setting turns out to
be the following.

\begin{condition} \label{main.cond}
  $X$ is a   scheme that is essentially of finite type over a field $k$,
 $x\subset X$ is a 0-dimensional closed subscheme, 
 $X$ satisfies Serre's condition  $S_2$ and 
 has pure dimension $\geq 3$ (that is,
every associated prime of $\o_X$ 
has the same dimension $\geq 3$). 
\end{condition}

\cite{MR492263}   defines a {\it local Picard functor} $\piclf(x, X)$
and proves that,
under the above assumptions, it is represented by the 
{\it local Picard scheme}
   $\picls(x,X) $; see  Definition \ref{loc.pic.func.defn}
and \cite[Thm.1.ii]{MR492263} for details. 
Frequently $\picls(x,X) $ has infinitely many connected components.
The identity component, usually
denoted by $\piclo(x,X) $, is  a commutative, connected algebraic group
of finite type but it need not be a linear algebraic group.
The union of those components of $\picls(x,X) $
that become torsion modulo
$\piclo(x,X) $ is denoted by $ \piclt(x,X)$.

The following theorem gives a unified answer to both questions.

\begin{thm}\label{main.thm.0}
Let $X$ and $Y$ be schemes that satisfy  Condition \ref{main.cond} and
 $\pi:Y\to X$  a finite morphism. 
Let $x\subset X$ be a $0$-dimensional closed subscheme,
  $y\subset Y$  the  preimage of $x$ and 
$r^X_Y: \picls(x,X)\to \picls(y,Y)$  the natural  pull-back map.

Assume that $\pi(Y)$ contains the support of a complete intersection
subscheme $x\subset Z\subset X$ of dimension $\geq 3$. 
Then  
$$
(r^X_Y)^{-1} \piclt(y,Y)=\piclt(x,X).
$$
\end{thm}

\begin{say}[Local N\'eron--Severi groups]  \label{loc.ns.intro}
One can restate  the above result
in terms of the {\it local N\'eron--Severi groups}
$$
\nsl(x,X):=\picls(x,X)/\piclo(x,X).
\eqno{(\ref{loc.ns.intro}.1)}
$$
The local N\'eron--Severi group is a finitely generated abelian group
if $\chr k=0$ but only the prime-to-$p$ part is known to be
finitely generated if $\chr k=p>0$; see 
Definition \ref{loc.NS.defn} for details.
Theorem \ref{main.thm.0}  is  equivalent to saying that
 the kernel of the natural pull-back map $\nsl(x,X){\longrightarrow}\nsl(y,Y)$
is torsion.

The global variant of this 
 was proved in \cite[p.305]{Kleiman66b}.

\end{say}

Theorem \ref{main.thm.0}
 implies that $\ker r^X_Y\subset \piclt(x,X)$ but for both Questions 
considered above 
we can say more. Our answers are optimal in characteristic $0$ but
in   characteristic $p>0$ I have not been able to
 exclude a possibly infinite,
discrete, $p^{\infty}$-torsion quotient.
The precise statements are the following.

\begin{thm}\label{pic-norm-loc.exc.thm1}
Assume that $(x,X)$ satisfies  Condition \ref{main.cond}, in particular
  $\dim X\geq 3$.
Let $\pi:Y\to X$ be a finite surjection
and $y\subset Y$ the  preimage of $x$. Then the kernel of 
the pull-back map
$\pi^*: \picls(x,X){\longrightarrow} \picls(y,Y)$ is
\begin{enumerate}
\item a linear algebraic group if $\chr k=0$ and
\item   the  extension of
a discrete $p^{\infty}$-torsion group by a
 linear algebraic group if $\chr k=p>0$.
\end{enumerate}
\end{thm}

I do not know any example where the kernel of $\pi^*$
 is not a linear algebraic group.

We have to be  careful with the scheme-theoretic formulation of
Question \ref{main.ques.2} since $D\subset X$ need not be $S_2$,
thus $\picls(x,D) $ need not exist.  Nonetheless 
for now I use the suggestive  notation
$r^X_D:\picls(x,X)
{\longrightarrow} \picls(x,D)$; see
Definition \ref{loc.pic.func.defn} for its precise meaning.

\begin{thm} \label{gl-gen.thm2}
Assume that $(x,X)$ satisfies  Condition \ref{main.cond} and  $\dim X\geq 4$.
Let  $x\subset D\subset X$ be an effective Cartier divisor. Then
the kernel of the restriction map
$r^X_D:\picls(x,X)
{\longrightarrow} \picls(x,D)$ is
\begin{enumerate}
\item a  unipotent  algebraic group if $\chr k=0$   or $X$ is normal and
\item  an extension of
a discrete $p^{\infty}$-torsion group by a
unipotent algebraic group if $\chr k=p>0$.
\end{enumerate}
\end{thm}

As before, 
I do not know any example where  the kernel  of $r^X_D $ is not unipotent.

\begin{rem}
In both cases the dimension restrictions are optimal.
For example, for the non-normal surface  $S:=(xyz=0)\subset \a^3$
with normalization $\pi:\bar S\to S$ we have 
$$
\ker\bigl[\picl(0,S)\stackrel{\pi^*}{\longrightarrow}
 \picls(\bar 0, \bar S) \bigr]\cong \z^3.
$$
For the ordinary 3-fold node  $X:=(x^2+y^2+z^2+t^2=0)\subset \a^4$
and $D:=(t=0)$ we have
$$
 \ker\bigl[\picls(x,X)\stackrel{r^X_D}{\longrightarrow}
 \picls(x, D) \bigr]\cong \z.
$$
There are also many instances where the kernels in 
(\ref{pic-norm-loc.exc.thm1}.1) and (\ref{gl-gen.thm2}.1)
are positive dimensional; see  Examples \ref{picl.from.res.cone.2},
 \ref{pic-norm-loc.exmp.2},
 \cite[Exmp.12]{k-gl1} and 
\cite[Exmps.1.31--35]{bha-dej}.

It is also important to assume in   (\ref{gl-gen.thm2}) that $D$ be a 
Cartier divisor. For instance, let $H\subset \p^n$ be a hyperplane
and $P\subset \p^n\setminus H$ a finite set of points. 
Let $(v,X)$ be a cone over the blow-up  $B_P\p^n$ with vertex $v$ and
$v\in D\subset X$ the cone over $H$. Then 
$$
 \ker\bigl[\picls(v,X)\stackrel{r^X_D}{\longrightarrow}
 \picls(v, D) \bigr]\cong \z^{|P|}.
$$
\end{rem}

If $\depth_xX\geq 3$ then 
 $\piclo(x,X) $ is trivial and we get the following results.

\begin{thm} \label{gl-gen.thm2.cor}
Notation and assumptions as in Theorems 
\ref{pic-norm-loc.exc.thm1}--\ref{gl-gen.thm2}.
Assume in addition that  $\depth_xX\geq 3$. Then
\begin{enumerate}
\item   $r^X_D: \picls(x,X){\longrightarrow} \picls(x,D)$ is injective,
\item   $\ker(\pi^*)$ is finite if $\chr k=0$  and
\item  the prime-to-$p$ part of 
 $\ker(\pi^*)$ is finite if $\chr k=p>0$.
\end{enumerate}
\end{thm}


\begin{say}[Numerical criteria for Cartier divisors]
\label{num.C.crit.ques.int.say}
Let  $T$ be an irreducible,  regular, 1-dimensional scheme and
 $f:X\to T$  a flat, projective morphism 
of relative dimension $n$. Assume for simplicity that $f$ has
  normal fibers.

Let $D$ be a divisor on $X$ such that $D_t:=D|_{X_t}$ is Cartier
for every $t\in T$. In general $D$ need not be Cartier.
For example, let 
$$
X=(x^2-y^2+z^2-t^2=0)\subset \a^3_{xyz}\times \a_t
\qtq{and} D=(x-y=z-t=0).
$$
  $D$ is Cartier, except at the origin,
where it is not even $\q$-Cartier.
However  $D_0$ is a line on a quadric cone, hence
$2D_0=(x-y=0)$ is Cartier. Thus $2D$  is 
 Cartier on every fiber but it is not Cartier.

We discuss several  criteria in Section
\ref{num.crit.sec}. The following is an easy-to-state special
case of Theorem \ref{num.C.crit.ques}.
\end{say}

 \begin{thm}\label{num.C.crit.ques.int.thm} Using the above notation,
assume in addition that $D_t$ is ample for every $t\in T$.
Then $D$ is a Cartier divisor on $X$ iff the self-intersection number
$(D_t^n)$ is independent of
$t\in T$.
\end{thm}

\begin{say}[The complex analytic case]  \label{anal.conjs}
Even if $k=\c$, the proofs of 
Theorems \ref{pic-norm-loc.exc.thm1}--\ref{gl-gen.thm2}
proceed by reduction to positive characteristic.
There are non-isolated complex analytic singularities
that do not lie on any algebraic variety. Our proof
does not apply to them, but the conclusions are hopefully 
valid. 

Over $\c$ one can use the first Chern class to realize
$\nsl(x,X) $ as a subgroup of  $H^2\bigl(\link(x,X), \z\bigr)$,
the second cohomology of the {\it link} of $(x,X)$, see 
(\ref{exp.seq.say}).
Thus it would be natural to try to prove the analytic case
by showing that  the kernels of
the corresponding maps between these cohomology groups are torsion.
This is, however, not true.

 Example \ref{pic-norm-loc.exmp.5}
 shows a semi-log canonical hypersurface singularity $(x,X)$
of dimension 3
whose normalization  $\pi:(\bar x,\bar X)\to (x,X)$
 is also a hypersurface singularity
and
$$
\ker \bigl[\pi^*: H^2\bigl(\link(x,X), \z\bigr){\longrightarrow}
H^2\bigl(\link(\bar x,\bar X), \z\bigr)\bigr]\cong \z^2.
$$
For Theorem \ref{gl-gen.thm2}, similar examples are given in 
\cite[Sec.5]{k-gl1}.
\end{say}

\section{Outline of the proofs}

The proofs of Theorems \ref{main.thm.0}--\ref{gl-gen.thm2.cor}
have 4 major components.
\begin{itemize}
\item The proof of  Theorem \ref{gl-gen.thm2} when $X$ is normal
and $\chr k>0$. This was done in \cite{bha-dej}.
\item The proof of  Theorem  \ref{pic-norm-loc.exc.thm1}
over finite fields. This is done in Section \ref{fin.field.sec}.
These two together imply  Theorem \ref{main.thm.0}
over finite fields.
\item A lifting argument that derives the general assertions
from  the finite field cases. This is done in Sections
\ref{using.relpic.sec}--\ref{qcoh.sec}.
\item Dealing with torsion elements in $\ker r^X_D$.
The nonexistence of prime-to-$p$ torsion 
follows from  \cite[XIII.2.1]{sga2}; see
Paragraph \ref{unip.last.step}. The $p^{\infty}$-torsion
is excluded in Section \ref{pf.of.7.sec} using a global argument.
\end{itemize}

The key to the results is to find a good answer to the following.

\begin{ques}\label{outline.ques.1}
 Let $(x,X)$ be a  local scheme 
and $L$ a line bundle on $X\setminus \{x\}$. How can one check if $L$ is in
$\piclt(x,X)$?
\end{ques}

In the global case, when  $Y$ is a proper scheme over a field $k$,
there is a simple numerical criterion:  a line bundle $L$ 
 is in $\pict(Y)$ iff $L$ has degree 0 on 
every reduced, irreducible curve
$C\subset Y$.

For the local Picard group there are no
proper curves to work with and I do not know  any similar
numerical criterion to identify
$\piclo(x,X)$ or $\piclt(x,X)$ in general.
(If $X$ is normal and it has a resolution of singularities $X'\to X$
then one can work on $X'$ and use the exceptional curves.
This was used in \cite{MR492263} to prove that  $\nsl(x,X)$
is finitely generated.)

Over a finite field the  $\f_q$-points of 
a finite type group scheme  form a  finite group,
which gives the following.

 \begin{claim}  Let   $(x,X)$ be a  local scheme  over a finite field $\f_q$
that satisfies  Condition \ref{main.cond}
and $L$ a  line bundle on $X\setminus \{x\}$.  Then
$L\in \piclt(x,X)$ iff $L$ is torsion.
\end{claim}

This leads to a somewhat roundabout way of proving that
a line bundle $L$ is in the connected component of a Picard group: 
\begin{itemize}
\item reduce everything to finite fields, 
\item check that we get  torsion line bundles and
\item lift back to the original setting.
\end{itemize}

The last step is the critical one; let us see it in more
detail for local schemes $(x,X)$ of finite type over a field $k$.
We may assume that $k$ is finitely generated over its prime field,
thus we can view $k$ as the function field of  an integral scheme $S$
that is of finite type over $\z$. By a suitable choice of $S$ we may even 
assume that we have
\begin{enumerate}
\item a scheme flat and of finite type $X_S\to S$ with a section
$\sigma:S\to X_S$ and
\item a line bundle $L_S$ on $X_S\setminus \sigma(S)$
such that
\item over the generic fiber we recover $(x,X)$ and $L$.
\end{enumerate}

The key technical result that we need is the following.

\begin{claim}\label{claim.13}
 $L\in \piclt(x,X)$ iff the set of closed points
$$
\bigl\{s\in S \mbox{ such that } L_S|_{X_s} \mbox{ is torsion}\bigr\}
\qtq{is Zariski dense in $S$.}
$$
\end{claim}

This is quite easy to prove if the local Picard groups
of the fibers $(x_s, X_s)$ are themselves fibers of a
``reasonable'' group scheme  $\picls_S(\sigma, X)$.
Even in the proper case, such relative Picard groups
exist only under some restrictions; see
 \cite[Chap.8]{blr} for a detailed discussion. In the local case,
the existence of such a group scheme  $\picls_S(\sigma, X)$ is not known.
We prove that, at least after replacing $S$ by a
dense open subset of $\red S$, there
 is a good enough approximation of $\piclo_S(\sigma, X)$
to make the rest of the proof work; see Section \ref{univ.sec}.

Constructing  $\piclo_S(\sigma, X)$
is relatively easy if 
the generic point  $s_g\in S$ has characteristic $0$
 since then $\piclo(x_g, X_g)$
is a {\em smooth} algebraic group. However, in positive  characteristic
we need to understand the obstruction theory of the Picard functor.
A delicate technical point is that the obstruction theory
is governed by  $H^2_x(X, \o_X)$ which is usually
infinite dimensional. A section of a coherent sheaf over $S$
vanishes at the generic point iff it vanishes at a Zariski dense set
of closed points, but this is no longer true for quasi-coherent sheaves;
see Example \ref{loc.free.but.bad}.
Thus the general theory does not exclude the possibility that
the dimension of $\piclo(x_s, X_s)$ jumps at every closed point.

Once this issue is settled, we complete the proofs
as follows.

Since $S$ is of finite type over $\z$, the residue fields
at closed points are all finite.

Over $\f_q$  we  prove  Theorem  \ref{pic-norm-loc.exc.thm1}
for  the normalization $\pi:\bar X\to X$  by factoring it as
$$
\pi:\bar X\stackrel{\pi_3}{\longrightarrow} X^{\rm wn}
\stackrel{\pi_2}{\longrightarrow}
 \red X\stackrel{\pi_1}{\longrightarrow} X
$$
where $X^{\rm wn} $ is the weak-normalization of $X$.
The nontrivial part is finiteness of the kernel  for
$\pi_3^*$. This is established
in   Proposition \ref{pic-loc-ff.prop}, using the
quotient theory of  \cite{k-q} and Seifert $\gm$-bundles
as in \cite[9.53]{kk-singbook}.
Besides proving Theorem \ref{pic-norm-loc.exc.thm1} over  $\f_q$,
this reduces the proof of  Theorem \ref{gl-gen.thm2} 
in positive characteristic
to the case when $X$ is normal. The latter was done in  \cite{bha-dej}.

For all the  theorems we use  Claim \ref{claim.13} to 
pass to arbitrary base fields $k$.

\begin{say}  \label{reduce.to.norm.say}
Here we show that Theorem  \ref{pic-norm-loc.exc.thm1}
is implied by the special case when $Y=\bar X$.
To see this, let $\sigma:\bar Y\to Y$ be the  normalization of $Y$.
It is enough to show that the kernel of the composite
$$
 \picls(x,X)\stackrel{\pi^*}{\longrightarrow} \picls(y,Y)
\stackrel{\sigma^*}{\longrightarrow}\picls(\bar y,\bar Y)
$$
is contained in  $\piclt(x,X) $.
The map $\sigma^*\circ \pi^*$ is also the composite of 
$$
\picls(x,X) \stackrel{\pi^*}{\longrightarrow} \picls(\bar x,\bar X)
\stackrel{\bar \pi^*}{\longrightarrow}\picls(\bar y,\bar Y).
$$
We already know that the kernel of $\pi^* $
is contained in $\piclt(x,X) $.

If $q:U\to V$ is a finite surjection between irreducible  normal schemes,
then for any line bundle $L$ on $V$ we have
$L^{\deg U/V}\cong \norm_{U/V}q^*L$.
Thus the kernel of $\bar \pi^*$ is  torsion and so 
the kernel of $\sigma^*\circ \pi^*=\bar \pi^*\circ \pi^*$ is  
also contained in  $\piclt(x,X) $. 
\end{say}

\begin{say}
\label{unip.last.step}
Notation and assumptions as in  Theorem \ref{gl-gen.thm2}.
Here we show that once we know that $\ker(r^X_D)\subset \piclt(x,X)$,
the remaining claims about unipotence follow.

Let $L\in \picl(x,X)$ be a nontrivial line bundle such that
$L|_D\cong \o_D$ and $L^m\cong \o_X$ for some $m>0$ not divisible by
$\chr k$. Choose $m$ to be the smallest with these properties.
Then $L$ determines a degree $m$, irreducible, cyclic cover
$\tau:\tilde X\to X$ that induces a trivial cover on $D\setminus \{x\}$.
Thus  $\tilde D:=\tau^{-1}(D)$ is a Cartier divisor such that
$ \tilde D\setminus \{\tilde x\}$ has $m>1$ connected components.
By \cite[XIII.2.1]{sga2} this is impossible.  

Thus  $\ker(r^X_D)$ is torsion free if $\chr k=0$
and contains only $p^{\infty}$-torsion if $\chr k=p>0$. 
Using Lemma \ref{no.torsion.i.gp}, this implies that
once $\ker(r^X_D)$ is known to be of finite type, it is
unipotent  if $\chr k=0$ and unipotent up-to
$p^{\infty}$-torsion if $\chr k=p>0$.
(Note that a finite $p^{\infty}$-torsion group is unipotent
if $\chr k=p>0$, so  $\ker(r^X_D)$ is  unipotent unless the
$p^{\infty}$-torsion part is not finitely generated.) 
\end{say}

\begin{lem}\label{no.torsion.i.gp} Let $G$ be a  locally of finite type, 
commutative algebraic group
over a field $k$. Assume that
$G$ is torsion free if $\chr k=0$
and contains only $p^{\infty}$-torsion if $\chr k=p>0$. 

Then $G^{\circ}$ is unipotent and $G/G^{\circ}$ is torsion free if $\chr k=0$
and contains only $p^{\infty}$-torsion if $\chr k=p>0$. 
\end{lem}

Proof. A connected commutative algebraic group $H$
 has a unique  connected subgroup
$0\subset H_l\subset H$ such that $H/H_l$ is an Abelian variety and
$H_l\cong H_m+H_u$ is a linear algebraic group where 
$H_m$ is multiplicative and $H_u$ is unipotent. 

The existence of $H_l$ is called Chevalley's theorem; the first published
proof is in \cite{MR0076427}, see \cite{MR3088271} for a modern treatment.
The decomposition $H_l\cong H_m+H_u$ is in most books on 
linear algebraic groups,
see for instance \cite[Thm.4.7]{borel}.

Note that both
 Abelian varieties and multiplicative  groups contain many torsion
 elements. These lift back to torsion elements in $H$
using the following elementary observation.

Let $H$ be a group and $K\subset H$ a central subgroup.
Assume that $\bar h\in H/K$ is $m$-torsion and $K$ is $m$-divisible.
Then $\bar h$ lifts to an  $m$-torsion element  $h\in H$. \qed

\section{Definition of local Picard groups}

The literature is very  inconsistent, there are at least
four variants  called the  local Picard group by some authors.

\begin{defn}[Local Picard group]\label{loc.pic.defn}
Let $X$ be a scheme and $x\in X$ a point. 
Assume for simplicity that $X$ is excellent and   $\depth_x\o_X\geq 2$. 

The  {\it local Picard group}
 $\picl(x, X) $ is a group whose 
 elements are $S_2$ sheaves $F$ on some neighborhood
$x\in U\subset X$ such that $F$ is locally free on $U\setminus \{x\}$.
Two such sheaves give the same element if they are isomorphic over some
 neighborhood of $x$. The product is given by the $S_2$-hull of the
tensor product.

One can also realize the local Picard group  as
$\pic(\spec \o_{x,X}\setminus\{x\})$ or as the direct limit of
 $\pic(U\setminus\{x\})$ as $U$ runs through all open 
Zariski neighborhoods of $x$. 
Usually it is necessary to take the limit.

If $X$ is normal and $X\setminus \{x\}$ is smooth 
(or locally factorial) then  $\picl(x, X) $
is isomorphic to  the divisor class group of $\o_{x,X}$.

In many contexts it is more natural to work with the 
{\it \'etale-local Picard group} $\picel(x, X) :=
\pic(\spec \o_{x,X}^h\setminus\{x\})$ where $\o_{x,X}^h $
is the Henselization of the local ring $\o_{x,X} $.
Alternatively, $\picel(x, X)$ is  the direct limit of
 $\picl(x', X')$ as $(x', X')$ runs through all \'etale 
 neighborhoods of $(x, X)$. Usually it is necessary to take the limit;
see (\ref{pic.not.glob.exmp}). 

Even for isolated singularities over $\c$, it is quite hard to
understand the relationship between  $\picl(x,X) $ and $\picel(x,X) $.
By \cite{MR1242007} there are many singularities such that
$\picl(x,X)=0 $ yet $\picel(x,X) $ is large.

Example \ref{pic.not.glob.exmp} shows that, 
for some rather simple singularities,
one always has $\picl(x,X)\neq \picel(x,X) $.

\end{defn}

\begin{defn}[Picard group of a local ring]\label{loc.pic.alg.defn}
 Let $(R,m)$ be a semilocal ring.
Assume for simplicity that $R$ is excellent and   $\depth_mR\geq 2$. 
We can define its local Picard group  $\picl(R,m)$
purely algebraically as follows. Its elements are
isomorphism classes of finite $R$-modules $M$ such that 
$\depth_mM\geq 2$ and 
$M_r$ is locally free of rank 1 over
$R_r$ for every non--zerodivisor $r\in m$. 
The product is given by the $S_2$-hull of the 
tensor product. 

If $(R,m)$ is the (semi)local ring of a point $x$ 
(or of a $0$-dimensional subscheme) on a scheme $X$
then $\picl(R,m)=\picl(x,X)$. 

In particular, $\picl(x,X)$ does not  depend
on our choice of the base scheme.

Assume next that $(x,X)$ is essentially of finite type over $(s,S)$. 
Then $k(x)$
is a finitely generated field extension of $k(s)$. Pick a
transcendence basis $\bar t_1,\dots, \bar t_m\in k(x)/k(s)$, lift these
back to $t_1,\dots, t_m\in \o_{x,X}$
and localize $\o_{s,S}[t_1,\dots, t_m]$ at the generic point of 
its intersection with $m_x$ to get  $(m_R, R)$. 
 We can now view $\o_{x,X}$ as 
an  essentially of finite type $R$-algebra.
The advantage is that now $\o_{x,X}/m_x$ is a
finite extension of $R/m_R$. Thus  $\o_{x,X}$ is
the localization of a finite type $R$-algebra at a closed point.

That is, in the study of  local Picard groups on
schemes of finite type, it is sufficient to work with
 closed points.  
\end{defn}

\begin{defn}[Local Picard functor and scheme] \label{loc.pic.func.defn}
\cite{MR492263}
Let $k$ be a field and  $(x, X)$ a local, Noetherian $k$-scheme.
Set $U:=\spec_X\o_{x,X}\setminus \{x\}$. 

For a local $k$-algebra $A$ consider the pre-sheaf
$$
A\mapsto  \pic\bigl((U\times_k\spec A)^h\bigr)
$$
where the superscript $h$ denotes the Henselization. 
Sheafifying in the \'etale topology gives the
{\it local Picard functor}  $\piclf(x, X)$.

Thus the local Picard functor works with objects
$\bigl( \pi:(\tilde x, \tilde X)\to (x,X), \tilde L\bigr)$
where $\pi$ is \'etale and $\tilde L $ is a line bundle on
$ \tilde X\setminus \tilde x $.

By \cite[Thm.1.ii]{MR492263}, 
if $\depth_xX\geq 2$ and 
$H^1(U, \o_U)\cong  H^2_x(X, \o_X)$ is finite dimensional, then 
the local Picard functor is represented by a  
  $k$-group scheme  $\picls(x,X) $ that is locally of finite type.
The  tangent space of  $\picls(x,X) $ at the identity is naturally isomorphic to
$H^2_x(X, \o_X)$.

If $H^2_x(X, \o_X)=0$ then $\picls(x,X)$ is essentially the same as
$\nsl (x,X)$. The interesting  case is when 
$0<\dim_k H^2_x(X, \o_X)<\infty$.  This holds if
\begin{enumerate}
\item[(i)]  $k(x)$ is finite over $k$,
\item[(ii)]  $X$ is $S_2$ and
\item[(iii)] $X$ is pure of dimension $\geq 3$.
\end{enumerate}
(These conditions are almost necessary, see  (\ref{push.finite.lem}).)

{\it Algebraic equivalence} is given by
\begin{enumerate}
\item[(a)] a connected $k$-scheme $T$ with two points $t_1, t_2$,
\item[(b)] an \'etale morphism  $\pi:Y\to X\times T$ such that the
injection $\{x\}\times T\into X\times T$ lifts to 
$\sigma: \{x\}\times T\into Y$ and
\item[(c)] a line bundle $L_Y$ on $Y\setminus \sigma\bigl(\{x\}\times T\bigr)$.
\end{enumerate}
Let a subscript ${\ }_i$ denote restriction to the fiber over $t_i$.
The two line bundles  
$$
\bigl(\pi_1: (y_1, Y_1)\to (x,X), L_1\bigr)\qtq{and}
\bigl(\pi_2: (y_2, Y_2)\to (x,X), L_2\bigr)
$$
are declared {\it algebraically equivalent.}
Note that  $L_1, L_2\in \picl(x,X)$ are  algebraically equivalent 
iff they are algebraically equivalent after some field extension $K\supset k$.
All line bundles algebraically equivalent
to the trivial bundle $\o_X$ form the identity component of $\picls(x,X) $,
denoted by $\piclo(x,X) $. As usual, 
$ \piclt(x,X)\subset  \picls(x,X)$
denotes the union of those components
that become torsion elements in
$\picls(x,X)/\piclo(x,X) $.

If  $H^2_x(X, \o_X)$ is infinite dimensional
then $\picls(x,X)$ does not exist, but, as long as $\depth_xX\geq 2$,
 the unit section of $\piclf(x, X)$ is represented by a
finitely presented closed immersion by   \cite[Thm.1.i]{MR492263}.

Assume now that $X$ is $S_2$,  $\dim X\geq 4$ and
 $x\in D\subset X$ is a Cartier divisor.
Set $\tilde D:=\spec_D j_*\o_{D\setminus x}$ and 
let $\tilde x\subset \tilde D$ be 
the preimage of $x$. Then $\depth_{\tilde x}\tilde D\geq 2$ 
hence  there is a subgroup scheme
${\mathbf K}(X|D)\subset \picls(x,X)$ representing those line bundles
 that become trivial when restricted to 
$D\setminus \{x\}\cong \tilde D\setminus \{\tilde x\}$.
This ${\mathbf K}(X|D)$ gives the precise definition of
$\ker\bigl[r^X_D:\picls(x,X)
{\longrightarrow} \picls(x,D)\bigr]$ used in Theorem \ref{gl-gen.thm2}.
\end{defn}

\begin{defn}[Local N\'eron-Severi group] \label{loc.NS.defn}
Let $k$ be a field and  $(x, X)$ a local, Noetherian $k$-scheme
such that  $\picls(x,X)$ exists. The
quotient
$$
\nsl(x,X):=
\picls(x,X)/\piclo(x,X)
$$
is called the {\it local N\'eron-Severi group.}
The analytic methods show that 
if $k$ has characteristic 0 then $\nsl(x,X) $
is a finitely generated abelian group, see (\ref{pic.H2.lem}). 
Combining the method of \cite{MR492263} with 
\cite{deJ-alt}
 implies that, even in characteristic $p>0$,
the local N\'eron-Severi group is finitely generated if $X$ is normal
and finitely generated modulo $p^{\infty}$-torsion in general.  
Our computations show that the   $p^{\infty}$-torsion part has bounded
exponent, but say nothing about finite generation.
I do not know any examples 
satisfying Condition \ref{main.cond}
where $\nsl(x,X) $ is not finitely generated.

The above is probably not a completely agreed-upon definition; another 
candidate is the usually smaller
$$
\picl(x,X)/\bigl(\picl(x,X)\cap \piclo(x,X)\bigr).
$$
\end{defn}

\begin{say}[Comparing $\pic$ and  $\mathbf{Pic}$]\label{comp.pic.pic.say}
Let  $(x, X)$ be a local, Noetherian scheme over a field $k$.
By definition there is a natural injection
$$
\picl(x,X)\into \picls(x,X)(k)
\eqno{(\ref{comp.pic.pic.say}.1)}
$$
which is usually not a surjection. 
If $k$ is algebraically closed, then, essentially by definition,
$$
\picel(x,X)\cong \picls(x,X)(k).
\eqno{(\ref{comp.pic.pic.say}.2)}
$$
If, in addition, $X$ is complete or henselian then
$$
\picl(x,X)\cong \picls(x,X)(k).
\eqno{(\ref{comp.pic.pic.say}.3)}
$$

\end{say}

\section{Examples of local Picard groups}

Here we discuss a series of examples of local Picard groups.
They show that the assumptions of 
Theorems \ref{pic-norm-loc.exc.thm1}  and \ref{gl-gen.thm2}
are essentially optimal.

\begin{exmp} \label{picl.from.res}
Let $(0,X)$ be an isolated singularity
with a resolution $p:Y\to X$ that is an isomorphism over $X\setminus\{0\}$.
Let $\{E_i:i\in I\}$ be the exceptional divisors. There is a natural
exact sequence
$$
0\to \tsum_i \z[E_i]\to \pic(Y)\to \picl(0,X) \to 0.
\eqno{(\ref{picl.from.res}.1)}
$$
(Left exactness follows from \cite[3.39]{km-book}, 
but there are many other ways to prove it.)

Assume in addition that there is an effective, exceptional divisor
$E$ such that $\o_Y(-E)$ is $p$-ample. 
(Any resolution obtained by repeatedly blowing up subschemes 
whose support lies over  $\{0\}$ has this property.) 
Let $L$ denote the line bundle  $\o_X(-E)|_E$; it is ample by assumption.

Let $nE$ denote the subscheme of $X$ defined by $\o_X(-nE)$.
There is an exact sequence
$$
0\to  L^n\stackrel{h\mapsto 1+h}{\longrightarrow} \o^*_{(n+1)E}
\to  \o^*_{nE} \to 1.
\eqno{(\ref{picl.from.res}.2)}
$$
Since $L$ is ample, $H^1(E, L^n)=H^2(E,L^n)$ for $n\gg 1$. Thus
the natural restriction maps
$\pic\bigl((n+1)E\bigr)\to \pic(nE)$
are isomorphisms for  $n\gg 1$. If $\o_{X}$ is  complete 
then 
$$
\pic(Y)\cong  \pic(nE) \qtq{for} n\gg 1.
\eqno{(\ref{picl.from.res}.3)}
$$
(This probably needs Grothendieck's existence theorem
\cite[III.5.1.4]{ega}.) 
Thus 
$$
\pic(0,X)=\pic(nE)/\tsum_i \z[E_i] \qtq{for} n\gg 1.
\eqno{(\ref{picl.from.res}.4)}
$$

For normal surfaces over a field of characteristic 0,
 (\ref{picl.from.res}.4) is used in  \cite{mumf-top} to  define the 
local Picard scheme given.
For $\dim X\geq 3$ it is an alternate way of constructing
 the local Picard scheme. Note, however, that the
scheme structures of the two sides can differ in
positive characteristic, as we see next.
\end{exmp}

\begin{exmp}[Cones] \label{picl.from.res.cone}
Let $W\subset \p^N_k$ be a smooth projective variety
and $(0,X)$ the affine cone over $W$;
see \cite[3.8]{kk-singbook} for our conventions on cones. 
By blowing up the vertex
we get a resolution $p:Y\to X$  where 
$$
X=\spec_k\tsum_{r\geq 0}H^0\bigl(W,\o_W(r)\bigr)\qtq{and} 
Y=\spec_W\tsum_{r\geq 0}\o_W(r).
$$
 Thus (\ref{picl.from.res}.3) suggests that
$$
\dim T_0\pics(Y)\stackrel{?}{=}H^1(Y, \o_Y)=
\tsum_{r\geq 0}^{\infty}H^1\bigl(W,\o_W(r)\bigr).
\eqno{(\ref{picl.from.res.cone}.1)}
$$
On the other hand,
$\dim T_0\picls(0,X) $ equals $H^1(U, \o_U)$ where $U=X\setminus \{0\}$.
Since $U=\spec_W\sum_{-\infty}^{\infty}\o_W(r)$,
we see that
$$
\dim T_0\picls(0,X)=\tsum_{-\infty}^{\infty} H^1\bigl(W,\o_W(r)\bigr).
\eqno{(\ref{picl.from.res.cone}.2)}
$$
Compared with (\ref{picl.from.res.cone}.1) we see 
extra summands for $r<0$. This is not a problem in characteristic 0
where these vanish by Kodaira's theorem. 
However, the two formulas can give different answers in 
positive characteristic.

Assume next that the characteristic is 0. Then
the above considerations give an exact sequence
$$
0\to  \tsum_{r> 0}H^1\bigl(W,\o_W(r)\bigr)\to
 \picls(0,X)\to\pics(W)/\z[\o_W(1)]\to 0.
\eqno{(\ref{picl.from.res.cone}.3)}
$$
Let  $H\subset W$ a smooth hyperplane section
and $D\subset X$  the cone over it. A similar computation gives 
an exact sequence
$$
0\to  \tsum_{r> 0}H^1\bigl(H,\o_H(r)\bigr)
\to \picls(0,D)\to \pics(H)/\z[\o_H(1)]\to 0.
\eqno{(\ref{picl.from.res.cone}.4)}
$$
Thus we see that the restriction  maps
$\picls(0,X)\to \picls(0,D)$ has a positive dimensional kernel
if the maps $H^1\bigl(W,\o_W(r)\bigr)\to H^1\bigl(H,\o_H(r)\bigr)$
are not all injective, cf.\ \cite[Exmps.1.31--35]{bha-dej}.
\end{exmp}

\begin{exmp}[Weighted cones] \label{picl.from.res.cone.2}
In (\ref{picl.from.res.cone}) let $(t=0)$ be an equation
of $H\subset W$. We can view $t$ as a map
$X\to \a^1$. After base-change to $t=s^m$ we get a new
singularity  $(0,X_m)$ containing the same $(0,D)$ as a
hyperplane section. Thinking of $X$ as
$\spec_k\tsum_{r\geq 0}t^rH^0\bigl(W,\o_W(r)\bigr)$, we now have
$$
X_m=\spec_k\tsum_{r\geq 0}s^rH^0\bigl(W,\o_W(\rdown{r/m})\bigr).
$$
Computing as in (\ref{picl.from.res.cone}) we obtain that
$$
H^2_0\bigl(X_m, \o_{X_m}\bigr)\cong
\tsum_{-\infty}^{\infty}H^1\bigl(W,\o_W(\rdown{r/m})\bigr)
\cong \bigoplus_1^m H^2_0\bigl(X, \o_{X}\bigr).
$$
Thus, if $\chr k=0$ then  
 $\dim \picls(0,X_m)=m\cdot \dim \picls(0,X)$ and so
$$
\dim\ker\bigl[\picls(0,X_m)\to \picls(0,D)\bigr]
$$
grows with $m$, save when $\picls(0,X)$ is 0-dimensional.

\end{exmp}

\begin{exmp} \label{pic.not.glob.exmp}
 Let  $S$ be a normal, projective  surface 
with a single isolated singularity at $(0, S)$. Assume that 
$\piclo(0, S)$ has no Abelian
subvarieties. We claim that   the image $\pic(S\setminus 0)\to \picl(0, S)$ 
is finitely generated.

To see this let $\pi:T\to S$ be a resolution of singularities. We can
factor $\pic(S\setminus 0)\to \picl(0, S)$  as 
$\pic(S\setminus 0)\to\pic(T)\to \picl(0, S)$. 
Since $\pic^0(T)$ is an Abelian variety, its image in
$\picl(0, S)$ is also an Abelian variety, thus trivial. Hence
$\pic^0(T)\to \piclo(0, S)$ is the constant map and so
$\pic(T)\to \picl(0, S)$ factors through  ${\rm NS}(T)\to \picl(0, S)$.

Here are some concrete equations with the above properties.

(\ref{pic.not.glob.exmp}.1) Cusps, for example  $\bigl(xyz+x^4+y^4+z^4=0\bigr)$
or $\bigl(z^2=x^2(x^2+y^2)+y^5\bigr)$. For these $\piclo(0, S)\cong \gm$.

(\ref{pic.not.glob.exmp}.2) Let $C$ be an irreducible curve with a
 single node whose
normalization $E:=\bar C$ is  a smooth elliptic curve.
There is an extension
$$
1\to \gm \to \pic^0(C)\to \pic^0(E)\cong E\to 0
$$
The extension is non-split, even more, 
$ \pic^0(C)$ does not contain any compact curves iff
the difference of the 2 preimages of the node is non-torsion on 
$\pic(E)$.
For example, the exceptional divisor of the minimal resolution of
$\bigl(z^2=x^2(x^4+cy^4)+y^7\bigr)$ is such a curve; the  non-torsion
condition holds for very general $c\in \c$.
\end{exmp}

In the non-normal case, several new phenomena occur.

\begin{exmp}[Non-normal surfaces] \label{pic-norm-loc.exmp.1}
We compute the local Picard group for three of the simplest
non-normal surfaces.  
Let $S:=(xy=0)\subset \a^3$ be the union of two planes in $\a^3$,
 $T:=(z^2=0)\subset \a^3$ the double plane in  $\a^3$
and $W:=(y^2-x^3=0)\subset \a^3$ the product of a cuspidal cubic with a line.
 $S,T$ and $W$ are all $S_2$.

Let $s\in S$ denote the origin; then $\depth_s\o_S=2$.
$S$ has two irreducible components  
$S_x=(x=0)\cong \a^2$ and $S_y=(y=0)\cong \a^2$.
The  normalization
$\bar S$ is the disjoint union $S_x\amalg S_y$. 
Thus $\picl(\bar s,\bar S)\cong 0$. We claim that
$$
\picl( s, S)\cong 
\ker\bigl[\picl( s, S)\stackrel{\pi^*}{\to} \picl(\bar s,\bar S)
\bigr]\cong \z.
\eqno{(\ref{pic-norm-loc.exmp.1}.1)}
$$
To see this let $L$ be an invertible sheaf
on $S\setminus s$. Both $L|_{S_x}$ and $L|_{S_y}$ are  trivial; let
$\sigma_x, \sigma_y$ be generating sections. Then $\sigma_x, \sigma_y$
restrict to global sections of $L$ on the punctured $z$-axis.
Their quotient is a regular function on the punctured $z$-axis;
its order of pole or zero  at the origin gives the isomorphism 
$\picl( s, S)\cong  \z$.

Next  let $t\in T$ denote the origin; then $\depth_t\o_T=2$.
The normalization of $T$ is  $\bar T=(z=0)\subset \a^3$.
There is an exact sequence
$$
0\to \o_{\bar T\setminus t}
\stackrel{e}{\to}  \o_{ T\setminus t}^*\to \o_{\bar T\setminus t}^*\to 1
$$
where $e(g)=1+zg$. Taking cohomology and using that $\picl(t, \bar T)=0$
gives the isomorphism
$$
\picl(t, T)\cong  H^1\bigl(\bar T\setminus t, \o_{\bar T\setminus t}\bigr).
\eqno{(\ref{pic-norm-loc.exmp.1}.2)}
$$
The latter  can be naturally identified with
$\sum_{m\geq 2}H^1\bigl(\p^1, \o_{\p^1}(-m)\bigr)$. Thus
$$
\begin{array}{rcl}
\picl(t, T)&\cong& 
\ker\bigl[\picl( t, T)\stackrel{\pi^*}{\to} \picl(\bar t,\bar T)
\bigr]\\
&\cong&  (\mbox{infinite dimensional vectorspace}).
\end{array}
\eqno{(\ref{pic-norm-loc.exmp.1}.3)}
$$

Let   $w\in W$ denote the origin.  We can write 
$W$ as $\spec k[t^2, t^3, z]$ with normalization
$\bar W\cong \spec k[t,z]$.  Let  $C\subset W$  denote the $z$-axis;
this is the singular locus of $W$.
There is an exact sequence
$$
1\to \o_W^*\to \o_{\bar W}^*\stackrel{d}{\longrightarrow} \o_C\to 0
$$
where  $d(f_0(z)+f_1(z)t+\cdots)=f_1(z)/f_0(z)$.
 Then we have the  exact sequence
$$
H^0\bigl(\bar W\setminus \{w\},  \o_{\bar W}^*\bigr)
\stackrel{d}{\longrightarrow}
H^0\bigl(C\setminus \{w\},  \o_C\bigr)\to
\pic(W\setminus \{w\})\to \pic(\bar W\setminus \{w\})=1.
$$
Note that $d$ factors as
$$
H^0\bigl(\bar W\setminus \{w\},  \o_{\bar W}^*\bigr)=
H^0\bigl(\bar W,  \o_{\bar W}^*\bigr)\stackrel{d}{\longrightarrow}
H^0\bigl(C,  \o_C\bigr)\to
H^0\bigl(C\setminus \{w\},  \o_C\bigr).
$$
This shows that  
$$
\begin{array}{rcl}
\picl(w\in W)&\cong & H^1_w(\o_C)\cong k[z,z^{-1}]/k[z]
\\
&\cong&  (\mbox{infinite dimensional vectorspace}).
\end{array}
\eqno{(\ref{pic-norm-loc.exmp.1}.4)}
$$
\end{exmp}

\begin{exmp} \label{pic-norm-loc.exmp.1.1}
As a  slight variation of  example (\ref{pic-norm-loc.exmp.1}), let
$X\subset \a^{2n+1}$ be the union of the linear spaces
$X_1=(x_1=\cdots=x_n=0)$ and $X_2=(x_{n+2}=\cdots=x_{2n+1}=0)$.
Here $X_1\cap X_2\cong \a^1$ and the normalization is
$\bar X=X_1\amalg X_2$. Here $X$ is not $S_2$ if $n\geq 2$ but
$\depth_x\o_X=2$ where $x\in X$ denotes the origin.

As in the previous example, we see that
$$
\picl( x, X)\cong 
\ker\bigl[\picl( x, X)\stackrel{\pi^*}{\to} \picl(\bar x,\bar X)
\bigr]\cong \z.
\eqno{(\ref{pic-norm-loc.exmp.1.1}.1)}
$$
\end{exmp}

\begin{exmp}[Demi-normal varieties]\label{pic-norm-loc.exmp.2} 
Fix a ground field $k$ an let $X_i$ be cones over the  Segre embedding of
$\p^1\times \p^n$ with vertices $v_i$. 
Note that the $X_i$ have rational singularities at the origin
and $\picl(v_i,X_i)\cong \z$.

Pick distinct points $p_j\in \p^1$ and
let $D_{ij}\subset X_i$ be the cone over  $\{p_j\}\times  \p^n$.
Note that $D_{ij}\cong \a^{n+1}$. 

The first example $Y_2$ is obtained by gluing $X_1$ to $X_2$ using the
natural  identifications
$D_{1j}\cong D_{2j}$ for $j=1,2$. 
Then $\bar Y_2=X_1\amalg X_2$ and so $\picl(\bar v, \bar Y_2)\cong \z^2$.
We claim that
$$
\begin{array}{l}
\picl( v, Y_2) \cong k^*+\z^2\qtq{and}\\
\ker\bigl[\picl( v, Y_2)\stackrel{\pi^*}{\to} 
\picl(\bar v,\bar Y_2)\bigr]\cong k^*.
\end{array}
\eqno{(\ref{pic-norm-loc.exmp.2}.1)}
$$
The extra  $k^*$ is obtained as follows. Set $D_i=D_{i1}+D_{i2}\subset X_i$. 
Take the trivial bundles
$\o_{X_i}$ and an isomorphism
$$
\phi^0: \o_{D_1\setminus v}\cong \o_{D_2\setminus v}
$$
to get an element of $ \picl(v,Y_2)$.
Two such isomorphisms give the same line bundle iff they differ
by multiplication by   sections of
$H^0\bigl(X_i\setminus v, \o^*_{X_i\setminus v}\bigr)$.
Since $\dim D_{ij}\geq 2$, the isomorphism $\phi^0$
extends to a pair of isomorphisms
$$
\phi_j: \o_{D_{1j}}\cong \o_{D_{2j}} \qtq{for $j=1,2$.}
$$
These  give isomorphisms
$$
\phi_j(v):k\cong k(v_1)\cong k(v_2)\cong k.
$$
The quotient $\phi_1(v)/\phi_2(v)$ gives a well defined element of $k^* $.

The second example $Y_3$ is obtained by gluing $X_1$ to $X_2$ using the
natural  identifications
$D_{1j}\cong D_{2j}$ for $j=1,2,3$. We claim that here the kernel is 
a direct product:
$$
\begin{array}{l}
\picl( v, Y_3) \cong k^{n+1}\times k^*+\z^2\qtq{and}\\
\ker\bigl[\picl( v, Y_3)\stackrel{\pi^*}{\to} 
\picl(\bar v,\bar Y_3)\bigr]\cong k^{n+1}\times (k^*)^2.
\end{array}
\eqno{(\ref{pic-norm-loc.exmp.2}.2)}
$$

To see this, set $D_i=D_{i1}+D_{i2}+D_{i3}\subset X_i$. 
 The trivial bundles $\o_{X_i}$ can be glued using a section of
$H^0\bigl(D_i\setminus v, \o^*_{D_i\setminus v}\bigr)$.
Two such sections give the same line bundle iff they differ
by multiplication by a  section of
$H^0\bigl(X_i\setminus v, \o^*_{X_i\setminus v}\bigr)$.
The previous considerations explain the $(k^*)^2 $
but now we also get new conditions from sections of the form
$1+(\mbox{linear functions})$.
On $D_i$, we have $n+1$ independent linear functions on each $D_{ij}$,
giving $3n+3$ independent linear functions  all together. 
On $X_i\subset \a^{2n+2}$ we only have
$2n+2$ independent linear functions. This accounts for the
$k^{n+1} $ part.

In principle we could get  more conditions by considering
$1+(\mbox{quadratic functions})$, but it is easy to see that this does
not happen for 3 or 4 copies of $D_{ij}$. For more copies, we get
conditions coming from higher degree polynomials.

\end{exmp}

\begin{exmp}[Pinch point]\label{pinch.exmp}
Let $k$ be any field and $S:=(x^2=y^2z)\subset \a^3$ the pinch point.
Its normalization is $\pi: \bar S:=\a^2_{uv}\to S$ given by
$(u,v)\mapsto (uv, u, v^2)$.  
The line $(x=z=0)$ generates $\picl(0,S)\cong \z/2$.

Note that $\pi$ is a homeomorphism
 if $\chr k=2$.  This suggests that in positive characteristic
there may not be a  perfect analog of the first Chern class as a mapping to
topological cohomology.

\end{exmp}

A more general version is the following.

\begin{exmp}\label{pic-norm-loc.exmp.3} Start with
$(\bar X,\bar D)\cong (\a^n,\a^{n-1} )$ and let $\bar x\in \a^{n-1}$
be the origin.
Let $\tau: \a^{n-1}\to\a^{n-1} $ be coordinate-wise multiplication
by a primitive $r$th root of unity. Assume that the characteristic does 
not divide $r$.
 
Construct $\pi:\bar X\to X$ by identifying the points in the $\tau$-orbits
 with each other. 
Thus $X$ is an affine variety, even CM.
For $r=2$  it has only double normal crossing
singularities along $D:=\pi(\bar D)$. 
We claim that
$$
\picl( x, X)=
\ker\bigl[\picl( x, X)\stackrel{\pi^*}{\to} 
\picl(\bar x,\bar X)\bigr]\cong \z/r.
$$
For any $r$th root of unity $\epsilon $,
the $\tau$-action on $\bar D$ can be lifted to $\o_{\bar D}$
 as  $g(x)\mapsto  \epsilon\cdot g\bigl(\tau^{-1}(x)\bigr)$.
Taking the quotient we get line bundles $L(\epsilon)$ over $X$
corresponding to  the $r$th roots of unity.

Choose coordinates such that $\bar D=(z_n=0)$.
A local trivialization of 
 $L(\epsilon)$ would correspond to an invertible function
$\phi(z_1,\dots, z_n)$ such that 
$\phi(\epsilon z_1,\dots, \epsilon z_{n-1},0)=
\epsilon \phi(z_1,\dots, z_{n-1},0)$. This would imply
$\phi(0,\dots, 0)=\epsilon \phi(0,\dots, 0)$, thus
$\phi(0,\dots, 0)=0$ if $\epsilon\neq 1$. Therefore $\phi$ is not invertible
near the origin if $n\geq 2$. 
\end{exmp}

\section{Analytic local Picard groups}\label{anal.sect}

\begin{defn}[Analytic local Picard groups]\label{loc.pic.anal.defn}
 Let  $X$ be a complex analytic space and $x\in X$ a point.
Assume for simplicity  that    $\depth_x\o_X\geq 2$. 

Let $W\subset X$ be the intersection 
of $X$ with a small (open) ball around $x$. The 
{\it analytic local Picard group} $\pical(x, X)$
can be defined as in (\ref{loc.pic.defn}) 
using (analytic) $S_2$ sheaves on $W$. 
By \cite{Artin69d}, if $X$ is an algebraic variety over $\c$ then
there is a natural isomorphism
$$
\picel(x, X)\cong \pical(x, X^{\rm an}).
\eqno{(\ref{loc.pic.anal.defn}.1)}
$$
By \cite{MR0245837}, if $X$ is $S_2$ and has pure dimension $\geq 3$ then
$$
\pica(W\setminus\{x\})\cong \pical(x, X).
\eqno{(\ref{loc.pic.anal.defn}.2)}
$$
Note that if $\dim X=2$ then $\pica(W\setminus\{x\})$ is 
infinite dimensional but  $\pical(x, W) $ is finite dimensional
if $X$ is normal. 
\end{defn}

\begin{say}[Exponential sequence]\label{exp.seq.say}
Let $U$ be a complex space. Then $\pic(U)\cong H^1(U, \o_U^*)$ and
the exponential sequence
$$
0\to \z_U\stackrel{2\pi i}{\longrightarrow} \o_U
\stackrel{\rm exp}{\longrightarrow} \o_U^*\to 1
$$
gives  an exact sequence
$$
H^1\bigl(U, \o_U\bigr)\to \pic(U) \stackrel{c_1}{\longrightarrow} 
H^2\bigl(U, \z\bigr).
$$
Let $X$ be  a complex space, $x\in X$ a  point and
set $U:=X\setminus \{x\}$. 
A piece of  the  local cohomology  exact sequence is
$$
 H^1\bigl(X, \o_X\bigr)\to 
H^1\bigl(U, \o_U\bigr)\to H^2_x\bigl(X, \o_X\bigr)\to H^2\bigl(X, \o_X\bigr).
$$
Thus if  $X$ is Stein  then we have an isomorphism
$$
H^1\bigl(U, \o_U\bigr)\cong H^2_x\bigl(X, \o_X\bigr).
$$
For local questions we should replace $X$ by a 
contractible open neighborhood  $x\in W\subset X$.
Then $W\setminus\{x\}$ is homotopy equivalent 
to $\link(x,X)$, which is the intersection of $X$ with a small
sphere centered at $x$. Thus 
$H^2\bigl(W\setminus\{x\}, \z\bigr)=H^2\bigl(\link(x,X), \z\bigr)$.

Combining with  \cite{MR0245837} and \cite{Artin69d} 
we obtain the following well known result.
\end{say}

\begin{prop} \label{pic.H2.lem} Let $(x, X$) be a 
$\c$-scheme of finite type
that is $S_2$ and has pure dimension $\geq 3$.
Then taking the first Chern class gives an exact sequence
$$
0\to \piclo(x,X)\to \picls(x,X)\stackrel{c_1}{\to}
H^2\bigl(\link(x,X), \z\bigr).
\qed
$$
\end{prop}

This suggests a
 topological way to proving Theorem \ref{pic-norm-loc.exc.thm1}
in characteristic 0:
we should prove that the kernel of the pull-back map
$$
 H^2\bigl(\link(x,X), \z\bigr)\to H^2\bigl(\link(\bar x,\bar X), \z\bigr)
$$
is torsion.
However, Examples \ref{pic-norm-loc.exmp.4}--\ref{pic-norm-loc.exmp.5} 
show  that  these maps can have non-torsion kernel, 
even for hypersurface singularities. We start with a projective example
and then we take a cone over it to get a local example.

\begin{exmp}[A singular K3 surface] \label{pic-norm-loc.exmp.4}
Let $g_4(x_0, x_1, x_2)$ be a general quartic form.
Then
$$
S:=\bigl(u^2=x_0^2g_4(x_0, x_1, x_2)\bigr)\subset \p^3(1,1,1,3)
$$
is a K3 surface with a double line  $L=(x_0=u=0)$.
Its normalization  $\pi:\bar S\to S$ is the smooth Del Pezzo surface of degree 2
$$
\bar S:=\bigl(v^2=g_4(x_0, x_1, x_2)\bigr)\subset \p^3(1,1,1,2).
$$
Thus $H^2(\bar S, \z)=\pic(\bar S)\cong \z^8$. 

The preimage of the line $L$ is a smooth elliptic curve
$E=(x_0=0)$. We claim that
$$
\ker\bigl[H^2(S, \z)\stackrel{\pi^*}{\to} H^2(\bar S, \z)\bigr]
\cong H^1(E, \z)\cong \z^2.
\eqno{(\ref{pic-norm-loc.exmp.4}.1)}
$$
To see this we start with the short exact sequence
$$
0\to \z_L\to \pi_*\z_E\to Q\to 0
$$
which shows that 
$$
H^1(L, Q)\cong  H^1(E, \z)\qtq{and} H^2(L, Q)=0.
\eqno{(\ref{pic-norm-loc.exmp.4}.2)}
$$
 The sheaf $Q$
also sits in the short exact sequence
$$
0\to \z_S\to \pi_*\z_{\bar S}\to Q\to 0
$$
and from this we get an exact sequence
$$
0\to H^1(L, Q) \to  H^2(S, \z)\to H^2(\bar S, \z) \to H^2(L, Q).
\eqno{(\ref{pic-norm-loc.exmp.4}.3)}
$$
Putting (\ref{pic-norm-loc.exmp.4}.2) and (\ref{pic-norm-loc.exmp.4}.3)
together gives (\ref{pic-norm-loc.exmp.4}.1).
\end{exmp}

\begin{exmp}\label{pic-norm-loc.exmp.5}
Using the previous notation, set
$$
X:=\bigl(u^2=x_0^2g_4(x_0, x_1, x_2)\bigr)\subset \a^4.
$$
Its normalization  $\pi:\bar X\to X$ is given as
$$
\bar X:=\bigl(v^2=g_4(x_0, x_1, x_2)\bigr)\subset \a^4.
$$
Set $U:=X\setminus\{0\}$ and $\bar U:=\bar X\setminus\{\bar 0\}$.
 Note that $U$ (resp.\ $\bar U$)
is a Seifert $\c^*$-bundle over $S$   (resp.\ $\bar S$).
This implies that
$$
H^2\bigl(\bar U, \z\bigr)\cong \z^7\qtq{and}
H^2\bigl(U, \z\bigr)\cong \z^9.
$$
(The Seifert $\c^*$-bundle structure gives these modulo torsion.
Since both $X,\bar X$ are hypersurface singularities, 
$U, \bar U$ are simply connected, thus the above $H^2$ are torsion free.) 
Thus $\pi^*$ is not injective on $H^2$ and 
$$
\z^2\cong H^1(E,\z)\cong \ker \bigl[H^2\bigl(U, \z\bigr)
\stackrel{\pi^*}{\to}H^2\bigl(\bar U, \z\bigr).
\eqno{(\ref{pic-norm-loc.exmp.5}.1)}
$$
By contrast, we claim that $\picl(0, X)=0$ for very general $g_4$.
To see this note first that $\q$-Cartier divisors $D$ on $S$ correspond to
 Cartier divisors $\bar D$ on $\bar S$ for which $\bar D|_E$ is invariant
under the Galois involution of  $E\to L$. 

Since $\bar S$ is a Del~Pezzo surface, the pair $(E,\bar S)$
is obtained from a pair  $(E,\p^2)$ by blowing up 7 points
$P_1,\dots, P_7\in E$.  Thus  the image of
$\pic(\bar S)\to \pic(E)$ is generated by  $\o_{\p^2}(1)|_E$ and 
the $\o_E(P_i)$. For very general choice of the $P_i$ these 
are independent in $\pic(E)$ and the only  Galois invariant
divisor classes are given by the multiples of 
$\o_{\p^2}(3)|_E\otimes \o_E(-P_1-\cdots-P_7)$. 
This is also the pull-back of the hyperplane class under the projection
$S\to \p^2_{\mathbf x}$. Therefore $\picl(0, X)=0$.

\end{exmp}

\section{Schemes over finite fields}\label{fin.field.sec}

We start with  the following strengthening
 of Theorem \ref{pic-norm-loc.exc.thm1}  over finite fields.

\begin{thm} \label{pic-nomrm-loc.frf.thm}
Let $k$ be a finite field and $X$ an excellent
$k$-scheme of pure dimension $\geq 3$ that is  $S_2$
(or at least topologically $S_2$, see Definition \ref{top-S2.defn}). 
Let $x\in X$ be a  point with finite residue field $k(x)$.
Let $\pi:\bar X\to X$ denote the normalization and $\bar x$ the preimage of $x$.
Then the kernel of the pull-back map
$$
\pi^* : \picl(x,X)\to \picl(\bar x,\bar X) \qtq{is torsion.}
$$
\end{thm}

Note that, by our Examples  
\ref{pic-norm-loc.exmp.1}--\ref{pic-norm-loc.exmp.2},  the
assumptions that $\dim X\geq 3$, $X$ be topologically $S_2$
and $k(x)$ be finite are all  necessary.

\medskip

Proof. We factor the normalization $\pi:\bar X\to X$ as
$$
\pi:\bar X\stackrel{\pi_3}{\longrightarrow} X^{\rm wn}
\stackrel{\pi_2}{\longrightarrow}
 \red X\stackrel{\pi_1}{\longrightarrow} X
$$
and show that for each step  the kernel of the pull-back map
on the local Picard groups is torsion.

For $\pi_1: \red X\to X$ this is done in Lemma \ref{finite.pi1.lem}
and for $\pi_2: X^{\rm wn}\to \red X$  in 
Corollary \ref{finite.pi2.lem}.
Both of these results are well known and hold in much greater generality.
The most delicate part is finiteness of the kernel  for
$$
\pi_3^*: \picl\bigl( x^{\rm wn},  X^{\rm wn}\bigr)\to 
\picl\bigl(\bar x, \bar X\bigr).
$$
 This is established
in   Proposition \ref{pic-loc-ff.prop}, using the
quotient theory of  \cite{k-q} and Seifert $\gm$-bundles
as in \cite[9.53]{kk-singbook}.\qed

\begin{lem} \label{finite.pi1.lem}
Let $k$ be a  field of positive characteristic and $X$ a 
Noetherian $k$-scheme.
Then  the kernel of the restriction map
$\pic(X)\to \pic(\red X)$ is killed by some
power of $\chr k$.
\end{lem}

Proof.  Let $J\subset \o_X$ be an ideal sheaf such that $J^2=0$.
Set $X_J:=\spec \o_X/J$. Note that $j\mapsto 1+j$ identifies
$J$ with the kernel of  $\o_X^*\to \o_{X_J}^*$. This gives
an exact sequence
$$
H^1(X,J)\to H^1(X, \o_X^*)\to H^1(X_J, \o_{X_J}^*).
$$
Since $H^1(X,J) $ is a  $k$-vector space,
it is killed by multiplication by $\chr k$.

Using this for the powers of the ideal sheaf defining
$\red X$ we see that 
the kernel of the restriction map $\pic(X)\to \pic(\red X)$ is 
killed by multiplication by a power of $\chr k$.\qed
\medskip

The above proof even shows the following stronger result.

\begin{cor} \label{finite.pi1.lem.cor}
Let $k$ be a  field of positive characteristic and $(x,X)$ a 
local Noetherian $k$-scheme that is $S_2$ and has pure dimension $\geq 3$.
Then  the kernel of the restriction map
$\picls(x,X)\to \picls(x,\red X)$ is a connected unipotent group scheme.\qed
\end{cor}

\begin{defn} \label{wn.defn} Let $X$ be a scheme.
A {\it partial weak normalization} is a 
finite morphism  $\pi:X'\to X$ such that
\begin{enumerate}
\item $X'$ is reduced,
\item  $p$ is a finite,  universal homeomorphism and 
\item the induced morphism $X'\to \red X$ is an isomorphism 
at all generic points.
\end{enumerate}
$X$ is called  {\it weakly normal} if $X$ is reduced and its sole
partial weak normalization is  the identity $X\cong X$.

If $X$ is  excellent, more generally, if the normalization 
$\bar X$ is finite over $X$,
then there is a unique maximal partial weak normalization
$X^{\rm wn}\to X$, called the {\it weak normalization} of $X$.
Note that $X^{\rm wn}$ is  weakly normal.

If $X$ has residue characteristic 0, then the
weak normalization agrees with the seminormalization.
See \cite[Sec.I.7.2]{rc-book} for details.

A weakly normal scheme is $S_2$ iff it is connected in codimension 1
as in \cite{MR0142547}. 
Therefore, if  $X$ is $S_2$ then so is its weak normalization.
 \end{defn}

For the proof of the next lemma, see   \cite[Sec.6]{k-quot} 
or  \cite[Prop.35]{k-q}.

\begin{lem}\label{fuh.frob.lem} 
Let $k$ be a  field  of positive characteristic 
and $g:Y\to X$ a finite morphism of Noetherian $k$-schemes.
The following are equivalent.
\begin{enumerate}
\item  $g$ is a universal homeomorphism.
\item  $Y$ is dominated by a Frobenius twist
$F_q:X_q\to Y\stackrel{g}{\to} X$
for some power $q$ of $\chr k$.\qed
\end{enumerate}
\end{lem}

\begin{cor} \label{finite-uh.lem}
Let $k$ be a   field  of positive characteristic 
and $g:Y\to X$ a finite, universal homeomorphism of Noetherian $k$-schemes.
Then  the kernel of the pull-back map
$g^*:\pic(X)\to \pic(Y)$ is killed by some
power of $\chr k$. 
\end{cor}

Proof. By Lemma \ref{fuh.frob.lem}, there is
a $q=(\chr k)^m$ such that $g$   factors as
$$
F_q:X_q\to Y\stackrel{g}{\to} X.
$$
Thus it is enough to prove that the kernel of the pull-back map
$F_q^*:\pic(X)\to \pic(X_q)$ is killed by $q$.

If a line bundle $L$ is given by an open cover $\{U_i\}$
and transition functions $\{\phi_{ij}\}$ then
$F_q^*L$ can be given by the open cover $\{F_q^{-1}U_i\}$
and transition functions $\{F_q^*\phi_{ij}\}$.

As an abstract  scheme $X_q$ is naturally isomorphic to $X$
and under this isomorphism  $F_q^*\phi_{ij}=\phi_{ij}^q$.
Therefore the Frobenius pull-back  
$F_q^*L$ is isomorphic to  $L^{\otimes q}$ under this isomorphism.
Thus the kernel of $F_q^*:\pic(X)\to \pic(X_q)$
is precisely the set of $q$-torsion elements.\qed
\medskip

Since the weak normalization of an excellent scheme 
is a finite, universal homeomorphism, as a special case
we get the following.

\begin{cor} \label{finite.pi2.lem}
Let $k$ be a   field  of positive characteristic 
and $X$ an excellent $k$-scheme  with weak-normalization  $\tau:X^{\rm wn}\to X$.
Then  the kernel of the pull-back map
$\pic(X)\to \pic\bigl(X^{\rm wn}\bigr)$ is killed by some
power of $\chr k$. \qed
\end{cor}

In contrast with (\ref{finite.pi1.lem.cor}),
 the kernel of $\tau^*$ need not be connected, as shown by
Example \ref{pinch.exmp}.

The last step of the proof uses the
concept of topologically $S_2$ schemes
and the theory of finite, set theoretic
equivalence relations developed in \cite{k-q}.

\begin{defn}[$S_2$ and topologically $S_2$]\label{top-S2.defn} 
Recall that a scheme $X$ is $S_2$
if a finite morphism $g:Y\to X$ is an isomorphism provided
\begin{enumerate}
\item there is a closed subset $Z\subset X$ of codimension $\geq 2$
such that $g$ is an isomorphism over $X\setminus Z$ and
\item $Y$ has no associated primes supported in $g^{-1}(Z)$.
\end{enumerate}

Similarly, a scheme $X$ is {\it topologically $S_2$}
if a finite morphism $g:Y\to X$ is a
finite, universal homeomorphism provided
\begin{enumerate}
\item[(1')] there is a closed subset $Z\subset X$ of codimension $\geq 2$
such that $g$ is a finite, universal homeomorphism  over $X\setminus Z$ and
\item[(2')]  $Y$ has no irreducible components supported in $g^{-1}(Z)$.
\end{enumerate}

It is not hard to see that a pure dimensional scheme 
$X$ is  topologically $S_2$ iff the following holds.
\begin{enumerate}\setcounter{enumi}{2}
\item[(3')]  Let $\pi:U\to X$ be an \'etale morphism from a 
connected scheme $U$ and $Z\subset U$ a closed subscheme of 
codimension $\geq 2$. Then $U\setminus Z$ is connected.
\end{enumerate}

These imply that 
 an $S_2$ scheme is topologically $S_2$,
 a weakly normal scheme is  $S_2$ iff it is topologically $S_2$
 and
 the weak normalization of an   $S_2$ scheme  is topologically $S_2$.
We also use that by \cite[XIII.2.1]{sga2} a Cartier divisor
in a  topologically $S_2$ scheme is also topologically $S_2$.
\end{defn}

\begin{say}[Set theoretic equivalence relations]\label{steq.say}
(For more details, see \cite{k-q}.)

Let $X$ be an excellent scheme, $\pi:\bar X\to X$ its normalization
and  $R$ the normalization of $\bar X\times_X\bar X$. 
Together with the coordinate projections
 $\sigma_1, \sigma_2:R\rightrightarrows \bar X$ 
we have a {\it  finite, set theoretic
equivalence relation.}

Let $Q\subset R$ be a closed subscheme such that
$\sigma_1|_Q, \sigma_2|_Q: Q\rightrightarrows \bar X$ 
is also a   set theoretic
equivalence relation. By \cite[Lem.1.7]{k-q}, the
geometric quotient  $\bar X/Q$ exists and the geometric fibers
of $\bar X\to \bar X/Q$ are precisely the $Q$-equivalence classes.

In particular, if $Q=R$ then we get that
$\bar X/R\to X$ is a finite, universal homeomorphism.
If $X$ is weakly normal then $\bar X/R= X$.
\end{say}

\medskip

Finally we study what happens
as we go from a weakly normal scheme to its normalization.

\begin{prop} \label{pic-loc-ff.prop.gen}
Let $k$ be a field of  characteristic $p>0$ and $X$ an excellent
$k$-scheme of pure dimension $\geq 3$ that is  weakly normal and $S_2$. 
Let $\pi:\bar X\to X$ denote the normalization.
Let $x\in X$ be a closed point with  residue field $k(x)$.
Then  there is a  $k$-torus $ {\mathbb T} $ and a linear representation
$$
I_x:\ker\bigl[\picls(x,X)\stackrel{\pi^*}{\longrightarrow} 
\picls(\bar x,\bar X)\bigr]\to {\mathbb T}
$$
whose kernel is $p^{\infty}$-torsion.
\end{prop}

Proof. Let $L$ be a line bundle on $X\setminus \{x\}$
such that  
$\pi^*L$ is trivial on $\bar X\setminus \{\bar x\}$.

We can view $L$ as a $\gm$-bundle over $ X\setminus \{ x\}$,
thus $\pi^*L$ is a trivial $\gm$-bundle over $\bar X\setminus \{\bar x\}$.
It extends to a trivial $\gm$-bundle $\bar L$ over $\bar X$.
Since $X$ is excellent, $\bar X$ is finite over $X$ hence
$\bar L$ is of finite type over $X$.

Let $R$ be the normalization of $\bar X\times_X\bar X$.
Then $\sigma_1, \sigma_2:R\rightrightarrows \bar X$ is a set-theoretic finite
equivalence relation. Since $X$ is  weakly normal,
  the geometric quotient $\bar X/R$ equals $X$;
see Paragraph \ref{steq.say}.


Let $r\subset R$ denote the union of the reduced preimages of $\bar x$.
Since $\pi^*L$ is pulled-back from $X$, we have isomorphisms
$$
\phi_{R\setminus r}:\sigma_1^*\bigl(\pi^*L\bigr)\cong 
\sigma_2^*\bigl(\pi^*L\bigr)\qtq{over $R\setminus r$.}
\eqno{(\ref{pic-loc-ff.prop.gen}.1)}
$$
By assumption $\pi^*L $ is trivial, thus
$\phi_{R\setminus r} $ is can be viewed as an  isomorphism of 
two trivial $\gm$-bundles on $R\setminus r$.

Let $R^{[1]}\subset R$ be the union of the 
 irreducible components of dimension $\geq \dim X-1$
and $R^{(1)}\subset R$   the equivalence relation 
generated by $R^{[1]}$.
We check at the end 
that $\bar X/R=\bar X/R^{(1)}$.

Since $R^{[1]}$  is normal and every  irreducible component
has dimension $\geq 2$,   $\phi_{R\setminus r} $ extends to
an isomorphism
$$
\phi_R^{[1]}:\sigma_1^*\bigl(\bar L\bigr)\cong 
\sigma_2^*\bigl(\bar L\bigr)\qtq{over $R^{[1]}$.}
\eqno{(\ref{pic-loc-ff.prop.gen}.2)}
$$

The isomorphisms (\ref{pic-loc-ff.prop.gen}.2) define a
 $\gm$-equivariant extension of the relation
$R^{(1)}$ to a finite relation on $\bar L$.
Let $\bar R_L$ denote the equivalence relation generated by this extension.
In general, such an extension is a pro-finite   equivalence relation.

In our case the only non-finiteness can occur
over $x\in X$. For   $\bar x_i\in \bar x$
let $ \gm(\bar x_i)$ denote the multiplicative group scheme
of the residue field $k(\bar x_i)$.
For a point $r_{\ell}\in R^{(1)}$ lying over $x$ 
set $\bar x_i:=\sigma_1(r_{\ell})$ and $\bar x_j:=\sigma_2(r_{\ell})$.
(We allow $i=j$.)
Then $\phi_R^{[1]} $ gives an isomorphism  
$$
\phi_{ij\ell}: \gm(\bar x_i)\cong \gm(\bar x_j)
\eqno{(\ref{pic-loc-ff.prop.gen}.3)}
$$
that is defined over $k(r_{\ell})$. 
 
It is enough to prove the Theorem after a finite field extension.
We can thus replace $k$ by the composite of
the above $k(\bar x_i)$ and $k(r_{\ell})$. 

Let $\Gamma_x$ denote the graph whose vertices are the points
$\bar x_i$ and to each $r_{\ell} $ we add an edge connecting
$\bar x_i$ and $\bar x_j$. Fixing a base point $\bar x_0$,
compositions of the above $\phi_{ij\ell} $ define a homomorphism
$$
I_x(L): H_1(\Gamma_x, \z)\to \gm(k).
\eqno{(\ref{pic-loc-ff.prop.gen}.4)}
$$
The construction is compatible with field extensions and
Henselization, thus we get a homomorphism
$$
I_x: \ker\bigl[\picls(x,X)\stackrel{\pi^*}{\longrightarrow} 
\picls(\bar x,\bar X)\bigr]\to
\Hom\bigl(H_1(\Gamma_x, \z),\gm\bigr).
\eqno{(\ref{pic-loc-ff.prop.gen}.5)}
$$
Assume next that $L$ is in the kernel of $I_x$.  
Then $\bar R_L$ is  a finite   equivalence relation.
More precisely, each  $\bar R_L$-equivalence class over
$x$ contains exactly 1 point of $\gm(\bar x_0)$.

Since $\bar L$ is of finite type over $X$,
by \cite[Thm.6]{k-q} the geometric quotient
$\bar L/\bar R_L$ exists and it is a Seifert $\gm$-bundle over
$\bar X/R^{(1)}$ by \cite[9.48]{kk-singbook}.
(The statement there assumes that two other conditions
(HN) and (HSN) are also satisfied. These are, however, used
only to ensure that  the geometric quotient
$\bar L/\bar R_L$ exists. In our case existence is assured by
\cite[Thm.6]{k-q}, the rest of the proof then works.)
By \cite[9.53]{kk-singbook},  a 
 power  $\bar L^{\otimes q}$ of $\bar L$ descends to a  $\gm$-bundle over
$\bar X/R^{(1)}$ for some $q=p^r$.

Finally, let $R^{(2)}\subset R$ denote the union of
all irreducible components of dimension $\leq \dim X-2$.
Let $Z\subset X$ be the image of $R^{(2)}$. 
Then $Z$ has codimension $\geq 2$ in $X$ and $R$ agrees
with $R^{(1)}$ over $X\setminus Z$.
Thus  $\bar X/R^{(1)}\to X$ is a finite, universal homeomorphism
over $X\setminus Z$. Since $X$ is  $S_2$,
this implies that $\bar X/R^{(1)}\to X$ is a finite, universal homeomorphism
over $X$. 
By construction, it is an isomorphism at all generic points.
Since $X$ is weakly normal, this implies that
$\bar X/R^{(1)}= X$. Thus $L^{\otimes m}$ extends to a $\gm$-bundle over $X$,
hence it is trivial in $\picl(x,X)$.
 \qed

\begin{say}\label{pic-loc-ff.prop.gen.say}
Note that we have proved that $I_x$ is a homomorphism of algebraic groups but
the proof did not establish that its image is a closed algebraic subgroup of
$ {\mathbb T} $. (For instance we did not exclude the possibility that
$\ker \pi^*\cong \z$ and $I_x$ is an injection.)  

However, once we know that $\ker \pi^*$ in 
 Theorem \ref{pic-norm-loc.exc.thm1} 
is of finite type, 
Proposition \ref{pic-loc-ff.prop.gen}  implies that
$\ker \pi^*$ is in fact linear.
\end{say}

\begin{cor} \label{pic-loc-ff.prop}
Let $k$ be a finite field and $X$ an excellent
$k$-scheme of pure dimension $\geq 3$ that is  weakly normal and $S_2$. 
Let $x\in X$ be a closed point.
Then the kernel of the pull-back map
$$
\pi^* : \picl(x,X)\to \picl(\bar x,\bar X) \qtq{is torsion.}
$$
\end{cor}

In contrast with the previous steps,  the order of the torsion
kernel  need not be a power of $\chr k$.
\medskip

Proof. Let $L$ be a line bundle on $X\setminus \{x\}$
such that  
$\pi^*L$ is trivial on $\bar X\setminus \{\bar x\}$.
In (\ref{pic-loc-ff.prop.gen}.4) we constructed a homomorphism
$$
I_x(L): H_1(\Gamma_x, \z)\to \gm(k).
$$
If  $k$ is finite, then $\gm(k) $ is a torsion group.
Thus $I_x(L^m)$ is the trivial homomorphism for some $m>0$.
By (\ref{pic-loc-ff.prop.gen}) this implies that
$L^m$ is  $p^{\infty}$-torsion, hence
$L$ is a torsion element of $\picl(x,X)$.
\qed

\medskip

We can now prove  the following form of
Theorem \ref{main.thm.0} over finite fields.

\begin{thm}\label{main.thm.0.ff} 
Let $k$ be a finite field, $X, Y$  excellent
$k$-schemes of pure dimension $\geq 3$ that are  topologically $S_2$ and
 $\pi:Y\to X$  a finite morphism. 
Let $x\subset X$ be a finite subscheme and
  $y\subset Y$  the  preimage of $x$.

Assume that $\pi(Y)$ contains the support of a complete intersection
subscheme $x\subset Z\subset X$ of dimension $\geq 3$. 

Then the kernel of the natural map
$\picl(x,X)\to \picl(y,Y)$ consists of torsion elements.
\end{thm}

Proof. We use induction on the dimension of $X$. 
If $\dim X\leq 3$ then in fact $X=Z$ and so $\pi$ is a finite surjection.
This case is treated in Theorem \ref{pic-nomrm-loc.frf.thm}
and Paragraph \ref{reduce.to.norm.say}.

Thus assume that $\dim X\geq 4$. If we replace $Y$ by some
$Y'\to Y$, the kernel only increases. Thus we may assume that
$\pi(Y)=\supp Z$,  $Y$ is normal and $Y\to X$ factors through
the normalization $\bar X\to X$.  Using Theorem \ref{pic-nomrm-loc.frf.thm}
 we may thus also
 assume that $X$ is normal.

Since $Z\subset X$ is a complete intersection, it is also
a complete intersection inside some Cartier divisor $D\subset X$.
As we noted in Definition \ref{top-S2.defn}, the divisor $D$ is
also topologically $S_2$ by \cite[XIII.2.1]{sga2}. 
Using \cite{bha-dej} we can thus replace $X$ by $D$ and finish by induction.
\qed

\section{Using the relative Picard group}
\label{using.relpic.sec}

Theorems   
\ref{main.thm.0}--\ref{gl-gen.thm2}
assert that certain line bundles are  contained in
$\piclt(x,X)$. We plan to prove such results by
first establishing the claim over finite fields and then
going back to arbitrary fields. 
As we already noted in Question \ref{outline.ques.1},
we need a method to decide when a 
 line bundle on $X\setminus \{x\}$ is contained in
$\piclt(x,X)$.

\begin{say}[General set-up]  Theorems \ref{main.thm.0},  
\ref{pic-norm-loc.exc.thm1} and \ref{gl-gen.thm2}
can be formulated as follows.

\medskip
{\it Step 0} (Starting point). 
Let $p:(y,Y)\to (x,X)$ be a finite morphism of local $k$-schemes
and $L\in \picl(x,X)$ such that $\pi^*L$ is in $ \piclt(y,Y)$.
By passing to a power of $L$, we may assume that 
$\pi^*L$ is in $ \piclo(y,Y)$, that is, $\pi^*L$ is algebraically
equivalent to the trivial bundle.

We would like to prove that, under suitable conditions,
$L\in \piclt(x,X)$. 

\medskip
{\it Step 1} (Spreading out).
There is an integral $\z$-scheme of finite type $S$  such that we have
the following.
\begin{enumerate}
\item[(i)] There are $S$-schemes of finite type  $Y_S\to S$ and
$X_S\to S$ with sections $\sigma_Y:S\to Y_S$ and $ \sigma_X:S\to X_S$.
\item[(ii)] There is a finite morphism 
$p_S:Y_S\to X_S$ such that $\sigma_X=p_S\circ \sigma_Y$
and $\red p_S^{-1}\bigl(\sigma_X(S)\bigr)=\sigma_Y(S)$.
\item[(iii)] There is a line bundle $L_S$ on $X_S\setminus \sigma_X(S)$
such that, for every $s\in S$ the restriction $p_S^*L_S|_{Y_s}$ is algebraically
equivalent to the trivial bundle on $Y_s$.
\item[(iv)] There is a map to the generic point  $\spec k\to S$ such that
$$
(p:Y\to X, L)\cong \spec k\times_S\bigl(p_S:Y_S\to X_S, L_S\bigr).
$$
\end{enumerate}
There are several known results that say that certain good
properties of the generic fiber are inherited by all fibers,
at least over a dense open subset;
see \cite[IV.12]{ega} for long lists.
For example, if the 
generic fiber is $S_2$ and has pure dimension $d$ then,
 possibly after shrinking $S$,  we may
assume that every fiber is $S_2$ and has pure dimension $d$.

\medskip
{\it Step 2} (Over finite fields).   For every closed point $s\in S$ we
have a finite morphism $p_s:(y_s,Y_s)\to (x_s,X_s)$ 
and $L_s\in \picl(x_s,X_s)$ such that $p_s^*L_s$ is algebraically
equivalent to the trivial bundle.

This is an instance of the original problem over 
the residue field $k(s)$. 
Since  $S$ is $\z$-scheme of finite type, these 
 residue fields are all finite. Thus
Theorem \ref{main.thm.0.ff}   implies that, under our assumptions,  
$L_s^{m_s}$ is trivial for some $m_s>0$.

If $m_s=m$ is independent of $s$ then it is reasonable to expect that
$L^m\cong \o_X$. However, this is usually not the case.
If $\piclo(x,X)$ has positive dimension and $\chr k=0$ then
there are non-torsion line bundles $L\in \piclo(x,X)$.
For these, the  $m_s$  are not even bounded.

\medskip
{\it Step 3} (Lifting; easy case).
If $\depth_{x_s}X_s\geq 3$ then \cite[XI.3.16]{sga2} implies that 
$L^{m_s}\cong \o_{X}$ and we are done. This happens precisely
when $\piclo(x,X)\cong \spec k$, thus   $\piclt(x,X)$ is identified
with the torsion subgroup of  $\nsl(x,X)$.

\medskip
{\it Step 4} (Lifting;  hard case).  
At this point we can forget about  $Y$  and $\pi:Y\to X$ in the original set-up.
Thus from now on we have only the following data. 
\begin{enumerate}
\item[(i)] A local $k$-scheme $(x,X)$ and $L\in \picl(x, X)$.
\item[(ii)] A finitely generated $\z$-subalgebra  $A\subset k$
and its spectrum $S$.
\item[(iii)] An $S$-scheme of finite type 
$X_S\to S$ with a section  $ \sigma_X:S\to X_S$.
We set $x_s:=\sigma_X(s) $ for $s\in S$. 
\item[(iv)] A line bundle $L_S$ on $X_S\setminus \sigma_X(S)$
such that $L_s\in \picl(x_s, X_s)$ is torsion for
 every closed point $s\in S$.
\item[(v)] An isomorphism
$\bigl(x,X, L\bigr)\cong \spec k\times_S\bigl(\sigma_X(S), X_S, L_S\bigr)$.
\end{enumerate}
The proofs of  Theorems \ref{main.thm.0}, 
\ref{pic-norm-loc.exc.thm1} and \ref{gl-gen.thm2}
will be completed by the next result. 
\end{say}

\begin{thm} \label{main.torsion.thm}
Let $S$ be an integral scheme whose closed points are dense in $S$.
Let $f:X\to S$ be a flat morphism whose fibers are $S_2$ and have
pure dimension $\geq 3$. Let  $\sigma:S\to X$
be a section and  $L$ a line bundle on $X\setminus\sigma(S)$.

Assume that, for a dense set of  closed point $s\in S$,
 there is a natural number
$m_s\in \n$ such that $L_s^{m_s}\cong \o_{U_s}$ where $U_s:=X_s\setminus\{x_s\}$. 

Then    $L_{k(S)}\in \piclt(x_{k(S)}, X_{k(S)})$.
\end{thm}

\begin{rem} It is possible that 
the assumption on  closed points being dense
is not necessary. By Proposition \ref{striv.then.pict.prop}
 this holds if the universal deformation in
Proposition \ref{prorep.prop} admits an algebraization.
I am able to prove only a weaker version of this: an algebraization of
$\piclo$ over the generic fiber gives an algebraization of the
universal deformation over an open subset of $S$.  This is why I need
to assume that closed points are dense in $S$; an assumption that always holds
in  our applications.

\end{rem}

\section{Formal deformation theory of $\picl$}

\begin{say}
\label{formal.def.th.say}
 Let  $(s, S)$ be a local scheme, $(y,Y)$ a local Henselian scheme and
$g:Y\to S$  a flat, affine  morphism.

Fix a line bundle $L$ on  $Y_s\setminus \{y\}$.
Let  $ (A,m_A)$ be a local Artin   $\o_S$-algebra 
with residue field $k=A/m_A$. 
Set $Y_A:=Y\times_S\spec A$ and   $U_A:=Y_A\setminus \{y\}$.

Let  $\defor_S(L,A)$ denote the set of isomorphism classes of
line bundles on $U_A$ whose restriction to
$U_k$ is isomorphic to (the pull-back of) $L$.

We check the conditions of \cite{MR0217093} for the pro-representability
 of  the functor $A\mapsto \defor_S(L,A)$.
Consider an extension 
$$
0\to M\to B\to A \to 0
\eqno{(\ref{formal.def.th.say}.1)}
$$ such that $M^2=0$.
Correspondingly there is  an exact sequence
$$
0\to M\otimes \o_{U_s}\stackrel{m\mapsto 1+m}{\longrightarrow}
\o_{U_B}^*\to \o_{U_A}^*\to 1.
\eqno{(\ref{formal.def.th.say}.2)}
$$
Pick any  $\gamma_A\in H^0\bigl(U_A, \o_{U_A}^*\bigr)$. 
If the fibers of $f$ are $S_2$ then  $\gamma_A$ extends to a section
$\gamma'_A\in H^0\bigl(Y_A, \o_{Y_A}\bigr)$. Since $f$ is affine,
we can lift $\gamma'_A$ to $\gamma'_B\in H^0\bigl(Y_B, \o_{Y_B}\bigr)$
and then restrict to  $\gamma_B\in H^0\bigl(U_B, \o_{U_B}^*\bigr)$.
Thus 
$$
H^0\bigl(U_B, \o_{U_B}^*\bigr)\to H^0\bigl(U_A, \o_{U_A}^*\bigr)
\qtq{is surjective.}
\eqno{(\ref{formal.def.th.say}.3)}
$$
Thus the following is a piece of 
the  long exact cohomology sequence  of (\ref{formal.def.th.say}.2)
$$
0\to M\otimes H^1\bigl(U_s, \o_{U_s}\bigr)\to  
H^1\bigl(U_B, \o_{U_B}^*\bigr)\to H^1\bigl(U_A, \o_{U_A}^*\bigr)
\stackrel{\obs}{\longrightarrow}
  H^2\bigl(U_s, \o_{U_s}\bigr)
$$
where $\obs$ is called the {\it obstruction map.}
If (\ref{formal.def.th.say}.1) 
splits then the liftings form a principal homogeneous space
under $M\otimes H^1\bigl(U_s, \o_{U_s}\bigr) $; the latter
 is  independent of $A$.
Furthermore, if $\dim X_s\geq 3$ then $H^1\bigl(U_s, \o_{U_s}\bigr) $ is
finite dimensional by (\ref{push.finite.lem}).

It remains to understand liftings to $D:=B\times_AC$ where 
 $B\onto A$ and $C\to A$ are maps of Artin  $\o_S$-algebras.
A line bundle $L_D$ on $U_D$ is determined by
a  line bundle $L_B$ on $U_A$, a  line bundle $L_C$ on $U_C$
plus an isomorphism $\phi: L_C|_{Y_A}\cong L_B|_{Y_A}$.
Any different  isomorphism  is given by
$\gamma_A\cdot \phi$ where  $\gamma_A\in H^0\bigl(U_A, \o_{U_A}^*\bigr)$. 
By (\ref{formal.def.th.say}.3)  one can lift $\gamma_A$ to 
$\gamma_B\in H^0\bigl(U_B, \o_{U_B}^*\bigr)$. 
Thus instead of changing $\phi$ by $\gamma_A$ we can change
$L_B$ by the isomorphism  $\gamma_B :L_B\to L_B$
to conclude that $L_D$ does not depend on $\phi$. That is
$$
\defor_S\bigl(L, B\times_AC\bigr)\cong
\defor_S\bigl(L, B\bigr)\times_{\defor_S(L, A)}\defor_S\bigl(L, C\bigr).
$$
We have thus proved the following.
\end{say}

\begin{prop}\label{prorep.prop}
 Let  $(s, S)$ be a local scheme, $(y,Y)$ a local Henselian scheme and
$g:Y\to S$  a flat  morphism  with $S_2$ fibers of dimension 
$\geq 3$. Let    $L$  be a line bundle  on  $Y_s\setminus \{y\}$.
Then $A\mapsto \defor_S(L,A)$ is pro-representable by a 
complete, local, Noetherian
   $\o_S$-algebra  $\deforb_S(L)$.\qed
\end{prop}

\begin{defn}\label{2.formal.defs.say}
 Let $S$ be a (non-local) scheme and  $f:X\to S$  a flat
 morphism  with $S_2$ fibers of pure dimension $\geq 3$.
Let $\sigma:S\to X$ be a section. For every $s\in S$ we
 are especially interested in  two deformation functors.

(\ref{2.formal.defs.say}.1) $\deforb_{(s,S)}(\o_{X_s})$  is obtained by applying
Proposition \ref{prorep.prop}
to $(s,S):=$ the localization of $S$ at $s$ 
and $(y,Y):=$ the Henselization of $X$ at  $\sigma(s)$.  
This is an $\o_{s,S}$-algebra
parametrizing  deformations 
(over Artin algebras over $\o_{s,S}$) of the trivial line bundle
$\o_{X_s} $ (pulled back to  $Y_s\setminus \{y\}$). 

(\ref{2.formal.defs.say}.2) $\deforb_{k(s)}(\o_{X_s})$ is obtained by applying
Proposition \ref{prorep.prop}
to $S:=\{s\}$ 
and $(y,Y):=$ the Henselization of the fiber $X_s$ at  $\sigma(s)$.  
This is a $k(s)$-algebra
parametrizing deformations 
(over Artin algebras over the residue field $k(s)$)
of the trivial line bundle
$\o_{X_s} $ (pulled back  to  $Y_s\setminus \{y\}$).

Note that 
$$
\deforb_{k(s)}(\o_{X_s})\cong \deforb_{(s,S)}(\o_{X_s})\otimes_Sk(s).
\eqno{(\ref{2.formal.defs.say}.3)}
$$
\end{defn}

\begin{defn}[Universal families]\label{vers.fam.defn}
 Let $S$ be a  scheme,  $f:X\to S$  a flat morphism  with $S_2$ fibers and
 $\sigma:S\to X$  a section. A {\it family} in $\piclf$ is given by
\begin{enumerate}
\item  a morphism $p:P\to S$,
\item an \'etale morphism  $g:Y\to X\times_SP$,
\item a lifting of the closed subscheme  
$(\sigma,1_P):P=S\times_SP\into X\times_SP$ to
$\sigma_P:P\to Y$ and
\item a line bundle  $L_Y$ on $Y\setminus \sigma_P(P)$.
\end{enumerate}
For a point $z\in P$ the corresponding line bundle 
on $Y_{p(z)}\setminus \sigma_P(z)$ is denoted by $L_z$.

These data give a {\it  finite type family} in $\piclf$
if $P$ is of finite type over $S$.

It is probably not necessary but for technical reasons it
is easier to assume that 
\begin{enumerate}\setcounter{enumi}{4}
\item $\red p:\red P\to \red S$ is smooth.
\end{enumerate}
This always holds in characteristic 0. In positive characteristic
it can be achieved after a purely inseparable, dominant
(but not necessarily finite) base change $S'\to S$
and by passing to an open subset $P'\subset P$. 

The above family is a  {\it deformation of the trivial bundle}
if, in addition, 
\begin{enumerate}\setcounter{enumi}{5}
\item  there is a section $\rho:S\to P$ such that the restriction of
$L_Y$ to $Y_{\rho}\setminus \sigma_P(P)$ is isomorphic to the
structure sheaf where  $ Y_{\rho}:=g^{-1}\bigl(\pi_2^{-1}(\rho(S))\bigr)$.
\end{enumerate}

For every $s\in S$ we can localize
at $s$ and get a deformation of $\o_{Y_s}$ as in (\ref{formal.def.th.say}). 

We say that the above family is {\it universal}
at a line bundle $L_s$ on $X_s\setminus \sigma(s)$
corresponding to a point $z\in P_s$ 
if the induced map
$$
\deforb_{(s,S)}\bigl(L_s\bigr)\to \widehat{\o}_{z,P} 
\qtq{is an isomorphism.}
\eqno{(\ref{vers.fam.defn}.7)}
$$
\end{defn}

\begin{prop} \label{striv.then.pict.prop}
Let $(s,S)$ be an integral,  local scheme with function field $K$.
Let  $f:X\to S$ be  a flat morphism
  with $S_2$ fibers,
 $\sigma:S\to X$  a section and $L$ a line bundle on $X\setminus \sigma(S)$.
Assume that
\begin{enumerate}
\item there is a  finite type family as in (\ref{vers.fam.defn}.1--6) that is
universal at $\o_{X_s}$ and
\item $L_s\cong \o_{X_s}$.
\end{enumerate}
Then  $L_K\in\piclt(x_K, X_K)$.
\end{prop}

Proof. For every Artin $\o_{X_s}$-algebra $A$, the line bundle 
$L$ defines  a deformation of $\o_{X_s}$ over $A$
since $L_s\cong \o_{X_s}$. This gives a formal deformation of 
$L_s\cong \o_{X_s}$ over  $\hat S$, the completion of $S$ at $s$.

By assumption 
$\deforb_{(s,S)}\bigl(\o_{X_s}\bigr)\cong \widehat{\o}_{\rho(s),P} $,
 thus  there is a section
$u:\widehat S\to P$ such that $(u^*L_Y)|_{X_A}$ is isomorphic to  $L|_{X_A}$ 
for every Artin $\o_{s,S}$-algebra $A$. 
Proposition \ref{form.triv.then.triv.lem} then implies that
$u^*L_Y\cong  L|_{\hat S}$.  Since $S$ is integral, 
the section $u$ factors through  $\red P$.


Now we use that $\red P\to  S$ is smooth (\ref{vers.fam.defn}.5).
Thus $\widehat{\o}_{\rho(s),\red P}$ is
 a formal power series ring over $\widehat{\o}_{s,S}$ and so
 each formal section of  $\red P\to  S$ can be 
given by an ideal
$\bigl(t_1-g_1,\dots, t_m-g_m\bigr)$ where $g_i\in \widehat{\o}_{s,S}$.
 Thus any two sections
are algebraically equivalent.
This implies that 
  $L_K$ and $\o_{X_K}$ are algebraically equivalent over
 $k(\hat S)$, hence also over $K$. \qed

\begin{prop} \label{form.triv.then.triv.lem}
Let $(s,S)$ be a local scheme with  maximal ideal $m$.
Let $f:X\to S$ be a scheme, flat over $S$ with $S_2$-fibers.
Let $X_n:=\spec_X \o_X/m^{n+1}\o_X$ be   the $n$th infinitesimal 
neighborhood of $X_0:=X_s$.
Let  $Z\subset X$ be a subscheme that is 
finite over $S$ and 
 $j:X\setminus Z\into X$ and $j_n:X_n\setminus Z_n\into X_n$
 the natural injections.
Let $L$ be a  locally free sheaf on $X\setminus Z$ and $L_n:=L|_{X_n\setminus Z_0}$.
 Assume that one of the following holds.
\begin{enumerate}
\item $(j_n)_*(L_n)$ is locally free for every $n\geq 0$.
\item $(j_0)_*(L_0)$ is locally free and $R^1(j_0)_*(L_0)=0$.
\end{enumerate}
Then $j_*L$ is  locally free in a neighborhood of $Z_0$.
\end{prop}

Proof. We may assume that $\o_S$ is $m$-adically complete  and, possibly after
passing to a smaller neighborhood of $Z_0$, we may assume that 
$f$ is affine and $(j_0)_*(L_0)\cong \o_X $.\
For every $n$  we have an exact sequence
$$
0\to (m_0^{n}/m_0^{n+1})\otimes L_0\to  L_n\to L_{n-1}\to 0.
$$
Pushing it forward   we get an exact sequence
$$
\begin{array}{l}
0\to  (m_0^{n}/m_0^{n+1})\otimes(j_0)_*(L_0)\to  (j_n)_*(L_n)
\stackrel{r_n}{\to} (j_{n-1})_*(L_{n-1})\to \\[1ex]
\hphantom{0} \to  (m_0^{n}/m_0^{n+1})\otimes R^1(j_0)_*(L_0).
\end{array}
$$
If  $(j_n)_*(L_n)$ is locally free then so is its restriction to $X_{n-1} $
 and $r_n$ gives a map of locally free sheaves
$$
\bar r_n: (j_n)_*(L_n)|_{X_{n-1}}\to (j_{n-1})_*(L_{n-1})
$$
that is an isomorphism on $X_{n-1}\setminus Z_{n-1}$. Since
$\depth_{Z_{n-1}}X_{n-1}\geq 2$, this implies that $\bar r_n$
is an isomorphism and so $r_n$ is surjective. 
The vanishing of $R^1(j_0)_*(L_0)$ also implies that  $r_n$ is surjective. 
Thus each $(j_n)_*(L_n) $ is locally free along $X_n$ and
the constant 1 section of  
$(j_0)_*(L_0)\cong \o_{X_0} $ lifts back to a 
nowhere zero global section of $\varprojlim  (j_n)_*(L_n) $.
Hence $\varprojlim  (j_n)_*(L_n)\cong \o_X$.

Furthermore, we have a natural map
$j_*L\to \varprojlim  (j_n)_*(L_n)\cong \o_X$
that is an isomorphism on $X\setminus Z$. Since 
$\depth_{Z} j_*L\geq 2$, this implies that $j_*L\cong \o_X$. \qed
\medskip

The examples below show that going from formal triviality
of deformations to actual triviality is not automatic.

\begin{exmp} Let $(e,E)\cong (e, E')$ be an elliptic curve. 
Set $X:=(E\setminus\{e\})\times E'$ and $p:X\to E'$ the second projection.
Let $\Delta\subset X$ be the diagonal and $L=\o_X(\Delta)$.

For $p\in E'\setminus\{e\}$ the line bundle $L|_{X_p}$ is a
nontrivial element of 
$$
\pic(X_p\setminus\{e\})\cong \pic(E\setminus\{e\})\cong \pico(E).
$$
but $L|_{X_e}$ is trivial.

For $m\im \n$ let $X_m\subset X$ denote the $m$th infinitesimal thickening of
the fiber $X_1:=X_e$. We have exact sequences
$$
H^1\bigl(X_1, \o_{X_1}\bigr) \to
H^1\bigl(X_{m+1}, \o^*_{X_{m+1}}\bigr) \to
H^1\bigl(X_{m}, \o^*_{X_{m}}\bigr) \to
H^2\bigl(X_1, \o_{X_1}\bigr).
$$
Since $X_1\cong E\setminus\{e\}$ is affine, this shows that
$$
\pic(X_m\setminus\{e\})\cong \pic(E\setminus\{e\})\cong \pico(E).
$$
Thus $L|_{X_m}$ is trivial for every $m$.
\end{exmp}

\begin{exmp} Consider the family of smooth, affine surfaces
$$
S:=\bigl(x^2+y^2=1+t^2z^2\bigr)\subset \a^3_{xyz}\times \a^1_t.
$$
Set $D:=(x-1=y-tz=0)$ and 
$ L:=\o_{S}(D)   $.

$S_t$ is a hyperboloid for $t\neq 0$, thus $\pic(S_t)\cong \z$
is generated by $L_t$. For $t=0$ we get a cylinder
and  $\pic(S_t)\cong \z/2$
is generated by $L_0$. As in the previous example, we see that
$L^2$ is trivial on all infinitesimal neighborhoods of $S_0$.
\end{exmp}

\section{Existence of universal families} 
\label{univ.sec}

\begin{say} \label{def.spread.say}
Assume that we have a  field $K$ that is finitely generated over 
its prime field and a local scheme $(x,X)$ of finite type
over $K$. Write  $P_K:=\piclo(x,X)$.

As in (\ref{vers.fam.defn}.1--4)   there is a universal  family 
$$
\bigl(Y_K\stackrel{g_K}{\to} X_K\times P_K {\to} P_K,  L_{Y_K}\bigr).
\eqno{(\ref{def.spread.say}.1)}
$$

Everything in (\ref{def.spread.say}.1) can be defined over a
finitely generated subring of $A\subset K$, thus there is
an integral scheme $S$ of finite type over $\spec \z$
such that (\ref{def.spread.say}.1) is the generic fiber of a family
$$
\bigl(Y_S\stackrel{g_S}{\to} X_S\times P_S {\to} P_S,  L_{Y_S}\bigr).
\eqno{(\ref{def.spread.say}.2)}
$$
\end{say}

The statement and the proof of the next result
closely follow \cite{Artin74c}. I go through the 
details for two reasons.  I always found  \cite{Artin74c} rather
concise and, more importantly,  not all the assumptions of  \cite{Artin74c} are
satisfied in our case. There are two main differences.
The automorphisms groups of our objects are all infinite dimensional, but,
as it was already observed in \cite{MR492263}, this does not seem to cause any
problems. A more difficult point is that the obstruction spaces
are also  infinite dimensional.

For the knowledgeable reader, Theorem \ref{free+comm.main.ptop}
is the only part  not contained in \cite{Artin74c, MR492263}.

\begin{thm}[Openness of universality]\label{vers.open.thm}
Using the above notation, there is a dense open subset $T\subset S$ such that
$$
\bigl(Y_T\stackrel{g_T}{\to} X_T\times P_T {\to} P_T,  L_{Y_T}\bigr)
$$
is everywhere universal (as in  Definition \ref{vers.fam.defn}).
\end{thm}

\begin{rem} One can imagine that
 $$
P_T \quad  ''   {=} '' \quad \piclo(\sigma_T, X_T),
$$
but we do not claim this. The main reason is that  there are
families of line bundles  $L$ over $X_S\setminus \sigma(S)$
such that $L_s\in \piclo(x_s, X_s)$ at the generic point but
not at some special points. 
Consider for example the family
$$
X_s:=\bigl(xy=uv(u+v+s)\bigr)\subset \c^4
\qtq{and}
D_s:=(x=u=0)+(x=v=0).
$$
For $s\neq 0$ we see that  $D_s\sim (x=0)$ is trivial in
$\picl(0,X_s)$. For $s=0$ we can use \cite[2.2.7]{k-etc} to see
that $\picl(0,X_0)\cong \z^2$ with $ (x=u=0)$ and $(x=v=0) $
as generators. Thus $D_0$ gives a non-torsion element in
 $\picl(0,X_0)=\nsl(0,X_0)$. 
I do not know whether such
 points $s\in S$ can be Zariski dense or not.

Thus $P_T$ should be viewed as an open neighborhood of the
zero section in the (possibly nonexistent) $\piclo(\sigma_T, X_T) $.
\end{rem}

Proof. 
By generic flatness we may assume that $P\to S$ is flat.

First we prove that 
$\deforb_{(s,S)}\bigl(\o_{Y_s}\bigr)\to \widehat{\o}_{\rho(s), P} $ 
is an isomorphism 
iff the map between the fibers over the points 
$\deforb_{k(s)}\bigl(\o_{Y_s}\bigr)\to \widehat{\o}_{\rho(s), P_s} $ is 
 an isomorphism.
This is completely general and follows from 
(\ref{etale.if.on.fibers.lem}). 

Consider the   tangent map
$$
t_{P_K/K}:T_{\rho(K), P_K}\to R^1j_K\o_{U_K}   
$$
 defined in (\ref{tnagnet.map.say}.5). By assumption
$t_{P_K/K} $ is an isomorphism.
We prove in (\ref{tnagnet.map.prop})  that the tangent maps
$$
t_{P_s/s}:T_{\rho(s), P_s}\to R^1j_s\o_{U_s}   
$$
are isomorphisms for all $s$ in a  Zariski open subset of $S$.

 If $P\to S$ is smooth then  a simple algebra lemma 
(\ref{etale.to.regular.lem}) shows that
$\deforb_{k(s)}\bigl(\o_{Y_s}\bigr)\to \hat{\o}_{\rho(s), P_s} $ is 
an isomorphism
and we are done. If $\piclo(x_K,X_K)$ is a smooth group scheme 
then, possibly after shrinking $S$, we may assume
that $P\to S$ is smooth.
A group scheme over a field of characteristic 0 is always smooth,
thus the proof of Theorem
\ref{vers.open.thm}, and hence also  the proofs of Theorems 
\ref{main.thm.0}--\ref{gl-gen.thm2}, 
are complete if $\chr k=0$.

 Otherwise we need a more detailed study of obstruction theory;
this is started in (\ref{obstr.map.1.say}).

\begin{say}[Tangent map]\label{tnagnet.map.say}
We continue with the notation of (\ref{def.spread.say}).
Let $I_{S,P}\subset \o_P$ be the ideal sheaf of $\rho(S)\subset P$.
We identify $S$ with $\rho(S)$ and set $R:=\spec_P \o_P/I_{S,P}^2$.
The ideal sheaf of $S\subset R$ is denoted by $I_S$. 
These data are encoded in  a diagram
$$
\begin{array}{cccll}
Y_S & \into & Y_R &  & \\
\downarrow && g\downarrow\uparrow\sigma  
& & \mbox{where $g$ is flat,}\\
S & \into & R & \qtq{where} & \mbox{$\sigma$ is a section and}\\
|| && \downarrow && I_S^2=0\\
S & = & S &&
\end{array}
\eqno{(\ref{tnagnet.map.say}.1)}
$$
Set $U_R:=Y_R\setminus \sigma(R)$ 
with natural injection $j:U_R\into Y_R$
and  $U_S:=Y_S\setminus \sigma(S)$.
Let $L_R$ be a line bundle on $U_R$ and $L_S:=L_R|_{U_S}$. 
There is an exact sequence
$$
0\to I_S\otimes_S L_S\to L_R\to L_S\to 0.
\eqno{(\ref{tnagnet.map.say}.2)}
$$
Pushing it forward by $j$ we get
$$
0\to I_S\otimes_S j_*L_S\to j_*L_R\to j_*L_S
\stackrel{\partial}{\to} I_S\otimes_S R^1j_*L_S.
\eqno{(\ref{tnagnet.map.say}.3)}
$$
Assume now that $L_S\cong \o_{U_S}$ and $Y_S\to S$ has $S_2$ fibers. Then
$j_*L_S\cong \o_{Y_S}$ and the exact sequence becomes
$$
0\to I_S\otimes_S \o_{Y_S}\to j_*L_R\to \o_{Y_S}
\stackrel{\partial}{\to} I_S\otimes_S R^1j_*\o_{Y_S}.
\eqno{(\ref{tnagnet.map.say}.4)}
$$
Here $\partial$ factors through  $\o_{Y_S}\to \o_Z\cong \o_S$
thus, if $I_S$ is locally free over $\o_S$,  $\partial$ is equivalent to
the {\it tangent map}
$$
t_{R/S}: \shom_S(I_S, \o_S)\to R^1j_*\o_{Y_S}.
\eqno{(\ref{tnagnet.map.say}.5)}
$$
(Note that $ \shom_S(I_S, \o_S) $ is isomorphic to the
relative tangent sheaf of $R/S$ restricted to $S\cong \sigma(S)$.)

If we started with a set-up as in (\ref{def.spread.say})
then
$R^1j_*\o_{Y_S}\cong R^1j_*\o_{X_S}$ and we get the following.
\end{say}

\begin{prop}  \label{tnagnet.map.prop}
Let  $\bigl(Y_S\stackrel{g_S}{\to} X_S\times P_S {\to} P_S, 
S \stackrel{\rho_S}{\to} P_S, L_{Y_S}\bigr)$
be a deformation as in (\ref{def.spread.say}.2).
Assume that $I_{S,P}/I_{S,P}^2$ is free over $S$
and  $R^1j_*\o_{X_S}$ is free and commutes with base change.
Then the tangent map
$$
t_{P/S}: \shom_S(I_{S,P}, \o_S)\to R^1j_*\o_{X_S} 
$$
has constant rank over a dense open subset of $S$.\qed
 \end{prop}

We have used the following commutative algebra lemmas.

\begin{lem}\label{etale.if.on.fibers.lem}
 Let $(m, S)$ be a local ring and
$\phi:(m_1, R_1)\to (m_2, R_2)$ a map of local $S$-algebras.
\begin{enumerate}
\item If $R_1$ is complete and $\bar\phi: R_1/mR_1\to R_2/mR_2$
is surjective  then $\phi$ is surjective.
\item If $\phi$ is surjective, $\bar\phi$ is an isomorphism
and $R_2$ is flat over $S$ then  $\phi$ is an isomorphism.
\end{enumerate}
\end{lem}

Proof. For every $r\geq 1$ we have a commutative diagram
$$
\begin{array}{ccc}
(m^r/m^{r+1})\otimes R_1/mR_1 & \onto & (m^r/m^{r+1})\otimes R_2/mR_2 \\
\downarrow  && \downarrow \\
m^rR_1/m^{r+1}R_1 & \stackrel{a_r}{\to} & m^rR_2/m^{r+1}R_2.
\end{array}
$$
The vertical arrows are surjective hence so is $a_r $.
By induction on $r$ we obtain that 
$$
\bar\phi_r: R_1/m^{r+1}R_1\to R_2/m^{r+1}R_2 \qtq{are surjective.}
$$
Since $R_1$ is complete, we can pass to the inverse limit to conclude
that
$$
R_1=\varprojlim  R_1/m^{r+1}R_1\to \varprojlim R_2/m^{r+1}R_2 \qtq{are surjective.}
$$
This factors through the injection $R_2\into \varprojlim R_2/m^{r+1}R_2$
thus $\phi:R_1\to R_2$ is surjective.

For part (2), let $J$ be the kernel of $\phi$. We have an exact sequence
$$
0\to J\to R_1\stackrel{\phi}{\to} R_2\to 0.
$$
Since $R_2$ is flat over $S$, tensoring with $S/m_SS$ is also exact, thus
we get $$
0\to J/m_SJ\to R_1/m_SR_1\stackrel{\bar\phi}{\to}  R_2/m_SR_2\to 0.
$$
This implies that $J/m_SJ=0$ thus $J=0$ by the Nakayama lemma.\qed

\begin{lem}\label{etale.to.regular.lem} Let $(R,m)$ be a complete local ring
and $(S,n)$ a regular, local $(R,m)$-algebra. 
Then $R= S$ iff the natural maps 
$R/m\to S/n$  and
$ m/m^2\to n/n^2$ are isomorphisms. 
\end{lem}

Proof. If $R= S$ then clearly $R/m\to S/n$  and
$ m/m^2\to n/n^2$ are isomorphisms. 

Conversely, assume that  $R/m\to S/n$  and
$ m/m^2\to n/n^2$ are isomorphisms. By induction on $r$ we see that
the natural maps 
$R/m^r\to S/n^r$  are surjective. Among local rings  $(A,m_A)$ with fixed 
embedding dimension $\dim_{A/m}m_A/m_A^2$, the length of
$S/n^r$ is the largest possible since $(S,n)$ is regular. Thus each 
$R/m^r\to S/n^r$ is an isomorphism. Since $(R, m)$ is complete
this implies that $R=S$. \qed

\begin{say}[Obstruction map]\label{obstr.map.1.say}
Let $S$ be a base scheme, $R\to S$ a flat scheme with a section
$\rho:S\to R$ with ideal sheaf $I_S\subset \o_R$. Let
$T\subset R$ be a subscheme with ideal sheaf $I_T$.
Assume that $\rho(S)\subset T$ and $I_SI_T=0$.

Let $g:Y_R\to R$ be a flat morphism with a section
$\sigma:R\to Y_R$. By restriction we get
$Y_S\to S$ and $Y_T\to T$.  These data are summarized in the following diagram.
$$
\begin{array}{cccccll}
Y_S & \into & Y_T & \into & Y_R& & \\
\downarrow &&\downarrow && g\downarrow\uparrow\sigma 
&&\mbox{where $g$ is flat}\\
S & \into & T  & \into & R&\qtq{where} & \mbox{$\sigma$ is a section and} \\
|| && \downarrow&& \downarrow & & I_SI_T=0\\
S & = & S & = & S &&
\end{array}
\eqno{(\ref{obstr.map.1.say}.1)}
$$
Finally write   $U_R:=Y_R\setminus\sigma(R)$, $U_T:=Y_T\setminus\sigma(T)$,
 $U_S:=Y_S\setminus\sigma(S)$ and
$L_T$ be a line bundle on $U_T$. We would like to understand when
$L_T$ extends to a  line bundle $L_R$ on $U_R$.

We have an exact sequence
$$
0\to I_T\otimes_S \o_{U_S}\stackrel{m\mapsto 1+m}{\longrightarrow}
\o^*_{U_R}\to \o^*_{U_T}\to 1.
\eqno{(\ref{obstr.map.1.say}.2)}
$$
This gives
$$
R^1j_*\o^*_{U_R}\to R^1j_*\o^*_{U_T}
\stackrel{\partial}{\longrightarrow} I_T\otimes_S R^2j_*\o_{U_S}.
\eqno{(\ref{obstr.map.1.say}.3)}
$$
The  line bundle $L_T$  corresponds to a section
$ \o_S\to R^1j_*\o^*_{U_T}$;
composing with $\partial $ gives
$$
[L_T]: \o_S\to I_T\otimes_S R^2j_*\o_{U_S}.
\eqno{(\ref{obstr.map.1.say}.4)}
$$
Thus $L_T$ extends to a line bundle  $L_R$ iff  $[L_T] =0$. 

If $I_T/I_SI_T$ is free over $S$ then  $[L_T]$ is
equivalent to a map, called the {\it obstruction,} 
$$
\obs(L_T, R):  \shom_S(I_T, \o_S)\to R^2j_*\o_{Y_S}.
\eqno{(\ref{obstr.map.1.say}.5)}
$$
and $L_T$ extends to a line bundle  $L_R$ iff
$\obs(L_T, R)=0$. 

For us the following two consequences are especially important.
\medskip

{\it Claim \ref{obstr.map.1.say}.6.}
Assume that  $S$ is integral, $I_T/I_SI_T$ is free over $S$ and
$R^2j_*\o_{X_S}$ is free and commutes with base change over $S$. 
\begin{enumerate}
\item[i)]  The obstruction map
$\obs(L_T, R)$
has constant rank and commutes with base change  
over a dense open subset of $S$.
\item[i)] If $\obs(L_T, R)$ is injective then $L_T$ can not be extended
over any subscheme $T\subsetneq T'\subset R$.
\qed
\end{enumerate}
\end{say}

\begin{say}[End of the proof of Theorem \ref{vers.open.thm}]
\label{vers.open.thm.pf}
We return to the setting of (\ref{def.spread.say}.1)
and, in addition, we choose a scheme $W_K\supset P_K$ such that
$W_K$ is smooth at the identity and has the same
tangent space as $P_K$. (Thus $W_K=P_K$ if $P_K$ is smooth.)

As before, 
everything is defined over a
finitely generated subring  $A\subset K$, thus there is
an integral scheme $S$ of finite type over $\spec \z$
such that (\ref{def.spread.say}.1) is the generic fiber of a family
$$
\bigl(Y_S\stackrel{g_S}{\to} X_S\times P_S {\to} P_S,  W_S\supset P_S, 
L_{Y_S}\bigr).
\eqno{(\ref{vers.open.thm.pf}.1)}
$$
In order to apply (\ref{obstr.map.1.say}.6), 
let $\widehat{\o}_{S,W}$ denote the completion of the structure sheaf
of $W_S$ along the section $\rho(S)$, $\hat I_S\subset \widehat{\o}_{S,W}$
the ideal sheaf of $\rho(S)$ and $\hat I_P\subset \widehat{\o}_{S,W}$
the ideal sheaf of $P_S$. Finally set
$T:=\spec_S  \widehat{\o}_{S,W}/\hat I_P$ and
$R:= \spec_S  \widehat{\o}_{S,W}/\hat I_S\hat I_P$.

Shrinking $S$ if necessary, the following conditions
can be satisfied:
\begin{enumerate}
\item[i)] $T, R$ are flat over $S$ and $I_T/I_SI_T$ is free over $S$.
\item[ii)] The tangent map
$t_{T/S}: \shom_S(I_{S}, \o_S)\to R^1j_*\o_{X_S}$
is an isomorphism and commutes with base change. (This follows from
(\ref{tnagnet.map.prop}).)
\item[iii)]  $R^2j_*\o_{Y_S}$ is free and commutes with base change. 
(This follows by applying  Theorem \ref{free+comm.main.ptop}
to $Y_S\to P_S$ and  the sheaf $F=\o_{Y_S}$.)
\end{enumerate}

Thus, by  Claim \ref{obstr.map.1.say}.6, the line bundle
$L_S|_{T}$ does not  extended
over any subscheme $T\subsetneq T'\subset R$.
Therefore  $P_S\to S$ is universal along the zero section
$\rho:S\to P_S$.  Since $P_S\to S$ is a group scheme, this implies
that it is everywhere universal.\qed
\end{say}

\section{Quasi-coherent higher direct images}
 \label{qcoh.sec}

A coherent sheaf over a reduced scheme
is  free over a dense open set, but 
there are  (even locally free) quasi-coherent sheaves that are not 
free over any dense open subset.  Even worse, a nonzero section may vanish
at every closed point.

\begin{exmp} \label{loc.free.but.bad}
Let $M\subset \q$ be the $\z$-submodule consisting of all
$m/n$ such that $n$ has no multiple prime factors. Let $\tilde M$ be the
corresponding quasi-coherent sheaf over $\spec \z$. Then $\tilde M$ is
locally free of rank 1, but it is not free. In fact, every global section of it
vanishes at all but finitely many points of $\spec \z$. 
\end{exmp}

Our aim is to  show that such bad behavior does not happen for certain
higher direct images.

 \begin{say}[Direct image functors] \label{dirm.im.defn}
 We study  functors ${\mathcal H}$
 with the following properties.
 
 Given  a morphism $f:X\to S$  and a quasi-coherent sheaf $F$ on $X$,
 ${\mathcal H}(F)$ is a quasi-coherent sheaf on $S$ and
 for every $p:T\to S$ there are base-change maps
 $$
 p^*{\mathcal H}(F)\to {\mathcal H}\bigl(p_X^*F\bigr).
\eqno{(\ref{dirm.im.defn}.1)}
 $$
  The best known examples are  ${\mathcal H}:=R^if_*$.
 
  We say that  ${\mathcal H}(F)$ {\it commutes with base change} if
 the base-change map (\ref{dirm.im.defn}.1) 
is an isomorphism   for every  $p:T\to S$.
We say that ${\mathcal H}$ is {\it free} 
 if  ${\mathcal H}(F)$  is a free quasi-coherent sheaf.

 We say that ${\mathcal H}(F)$ is 
{\it generically free and commutes with base change}  if there is
 a nonempty  open set $S^0\subset S$ such that 
  ${\mathcal H}(F^0)$ is free and commutes with base change
where   $F^0$ denotes  the restriction of
  $F$ to $X^0:=f^{-1}(S^0)$. 
   \end{say}

The Cohomology and Base change theorem for proper morphisms
implies that if $S$ is reduced, $f:X\to S$ is proper and
$F$ is a coherent sheaf on $X$ then $R^if_*F$ is
generically free and commutes with base change for every $i$.
  In our applications we are especially 
interested in some cases where ${\mathcal H}(F)$ is
  only quasi-coherent, even though $F$ itself is coherent.  
Let us start with an example where some higher direct image is either
not generically free or does not commute with base change.

\begin{exmp} Let $E$ be an elliptic curve and  $L$  the
Poincar\'e bundle on $E\times E$. Set $X:=\spec_{E\times E}\sum_{i\geq 0}L^i$
and $g:X\to E\times E\to E$ a projection.
Note that $g$ is smooth, has fiber dimension 2 but it is
neither affine nor proper.

If $e\in E$ is not a torsion  point  then
$H^0(X_e, \o_{X_e})=H^1(X_e, \o_{X_e})=0$. If $e\in E$ is a torsion  point
of order exactly $m$ then there are natural identifications
$$
H^0(X_e, \o_{X_e})\cong H^1(X_e, \o_{X_e})\cong \tsum_{i\in m\n}\ k(e).
$$
Thus we see that $g_*\o_X=0$ is free but it does
nor commute with base change over any open set. By contrast,
$R^1g_*\o_X$ is a sum of skyscraper sheaves of infinite rank
supported at the torsion points. Thus it is not
generically free but it does commute with base change.
\end{exmp}

 The following are some basic examples that we use.

\begin{thm}[Generic freeness]\label{gen.flat.thm}
Let $f:X\to S$  be an affine morphism of finite type, $S$ reduced and   $F$  a
 coherent sheaf on $X$.
Then $f_*F$ is generically free and commutes with base change.
\end{thm}

Proof. Frequently this is stated as ``Generic flatness:''
 under the above assumptions,
$F$ is flat over a dense, open subscheme $S^0\subset S$;
 see, for example \cite[Lect.8]{mumf66}. 
It is stated as ``Generic freeness'' in \cite[Sec.14.2]{eis-ca}, but without
the commutation with base change.
However, both proofs show the stronger forms.
\qed

\begin{exmp} \label{free+comms.mooth.lem}
Let $X=\a^n_S\cong \spec_S\o_S[x_1,\dots, x_n]$, 
$Z\subset X$ the zero section
and  $j:X\setminus Z\into X$  the natural injection.
Then $f_*\circ R^ij_*\o_{X\setminus Z}$ is free and commutes with base change.

Note that $f_*\circ j_*\o_{X\setminus Z}=f_*\o_X$ if $n\geq 2$
and, for $i\geq 1$, the sheaf  $R^ij_*\o_{X\setminus Z} $ 
is supported on $Z\cong S$, thus $f_*$ is an isomorphism.
Furthermore, the only nonzero case is
  $R^{n-1}j_*\o_{X\setminus Z} $ which can be identified with
the quasi-coherent sheaf freely generated by
$$
\frac1{x_1^{a_1}\dots x_n^{a_n}}\qtq{for all}  a_1\geq 1,\dots, a_n\geq 1.
$$
(This is equivalent to the computation of  the groups 
$H^{n-1}\bigl(\p^{n-1}, \o_{\p^{n-1}}(m)\bigr)$ done in \cite[III.5.1]{hartsh},
see especially the 4th displayed formula on p.226.)
For more explicit references see  \cite[p.685]{eis-ca}.

More generally, if $f:X\to S$ is smooth along a section
 $Z\subset X$ and $j:X\setminus Z\into X$  is the natural injection
then $f_*\circ R^ij_*\o_{X\setminus Z}$ is generically free for $i\geq 1$. 
\end{exmp}

We use the following basic coherence theorem.

\begin{thm}\label{push.finite.lem}\cite[VIII.2.3]{sga2}
Let $X$ be an excellent scheme, $Z\subset X$ a closed subscheme,
$U:=X\setminus Z$ and $j:U\into X$  the  open
embedding. 
Assume in addition that  $X$ is locally embeddable into a
regular scheme. For a coherent sheaf $G$ on $U$  and $n\in \n$
the following are
equivalent.
\begin{enumerate}
\item $R^ij_*G$ is coherent for $i<n$.
\item $\depth_uG\geq n$ for every point $u\in U$
such that $\codim_{\bar u}(Z\cap \bar u)=1$.\qed
\end{enumerate}
\end{thm}

We need some elementary properties of direct image functors.

  \begin{lem} \label{free+comm.lem} Let 
$0\to {\mathcal H}_1(F)\to {\mathcal H}_2(F)\to {\mathcal H}_3(F)\to 0 $
be an exact sequence 
of direct image functors.
\begin{enumerate}
 \item If 
  and ${\mathcal H}_1(F),{\mathcal H}_3(F)$ are free 
then so is ${\mathcal H}_2(F)$.
\item If ${\mathcal H}_2(F)$ is free,
  ${\mathcal H}_1(F)$ is coherent and $S$ is reduced 
 then  $ {\mathcal H}_3(F)$ is generically  free.
\item If 
  and ${\mathcal H}_1(F),{\mathcal H}_3(F)$ commute with base change
then so does ${\mathcal H}_2(F)$.\qed
     \end{enumerate}
  \end{lem}
 
The following is the main new result of this section.
For our applications we need the case $F=\o_X$ and $r=1,2$
but the proof of the general version is similar.

 \begin{thm} \label{free+comm.main.ptop}   Let   $S$ be a reduced
scheme,  
$f:X\to S$   an affine morphism of finite type and 
$F$   a coherent sheaf on $X$.
Let $Z\subset X$ be a section,   $U:=X\setminus Z$ 
and $j:U\into X$  the natural injection.
Assume that  $F|_U$   is flat  over $S$
and each $F|_{U_s}$ is pure  dimensional and $S_r$. 

Then $f_*\circ R^rj_*(F|_U)$ is generically free and commutes with
 base change.\footnote{K.~Smith pointed out  that 
the assumptions
can be considerably weakened; see  \cite{ks-future}.} 
\end{thm}

Note that if each $F_s$ is $S_{r+1}$ then $f_*\circ R^rj_*(F|_U)$ is
coherent by (\ref{push.finite.lem}) but otherwise
$f_*\circ R^rj_*(F|_U)$ almost always has infinite rank.
\medskip

 Proof.  
By induction on $r$, starting with $r=0$. 
Here we do not assume that $F_s$ is  pure dimensional.
If $\dim F_s\leq 1$ then $\supp F$ is affine and
$f_*\circ j_*(F|_U)\cong (f|_U)_*(F|_U)$
 is generically free and commutes with base change
by (\ref{gen.flat.thm}).

If every associated prime of $F_s$ has dimension
$\geq 2$ then $j_*(F|_U)$ is coherent by (\ref{push.finite.lem}).
We claim that it 
generically  commutes with base change. A sheaf $G$ on $X$ such that
$G|_U\cong F$ equals $j_*(F|_U)$ iff it has no associated primes
supported on $Z$ and $\depth_ZG\geq 2$. Both of these are open conditions
in flat families by  \cite[IV.12.1.6]{ega}.
 Since  $\supp j_*(F|_U)$ is affine over $S$, we are done as above.

Even if $F_s$ is not pure dimensional, we have an exact sequence
$$
0\to F^{\leq 1}\to F|_U\to Q\to 0
$$
where $F^{\leq 1}$ is the largest subsheaf with $\leq 1$-dimensional support;
thus every associated prime of $Q_s$ has dimension
$\geq 2$.
By pushing forward, we get
$$
0\to (f|_U)_*F^{\leq 1}\to f_*\circ j_*(F|_U)\to f_*\circ j_*Q\to R^1j_*F^{\leq 1}=0.
$$
Thus $f_*\circ j_*(F|_U)$ is is generically free and commutes with base change
by (\ref{free+comm.lem}.1).
 
Thus assume that $r\geq 1$.
In this case   $R^rj_*(F|_U) $ 
is supported on $Z\cong S$, thus $f_*$ is an isomorphism
and we will drop it from the notation.

Let $n$ be the fiber dimension of $\supp F\to S$. If $n\leq r$ then
$R^rj_*(F|_U)=0$ and we are done. Thus assume from now on that
$n\geq r+1$. We can also replace $X$ by  $\supp F$.

After possibly shrinking $S$, we can 
use Noether normalization to obtain a finite morphism
$g:X\to \a^n_S$ such that $Z$ is mapped isomorphically to the 0-section
$Z_0\subset \a^n_S$. Let $j_0:U_0:=\a^n_S\setminus Z_0\into \a^n_S$
be the natural injection.

After possibly further shrinking $S$ we may also assume that
$g^{-1}(Z_0)$ is the disjoint union of $Z$ and of another closed
subscheme $Z'$. Set $U':=X\setminus Z'$ with natural injection
$ j':U'\into X$.  Then 
$$
R^r(j_0)_*\bigl(g_*(F|_{U_0})\bigr)
=\bigl(g_* R^rj_*(F|_U)\bigr)+\bigl(g_* R^rj'_*(F|_{U'})\bigr).
$$
Thus if $R^r(j_0)_*\bigl(g_*(F|_{U_0})\bigr) $ is
 generically free and commutes with base change
then the same holds for its direct summands.
 Since $ R^rj_*(F|_U)$ is supported on $Z$ and
$g|_Z:Z\to Z_0$ is an isomorphism, $g_* R^rj_*(F|_U)$
is naturally isomorphic to $R^rj_*(F|_U)$. 

Thus we have reduce everything to 
the case when $f:X\cong \a^n_S\to S$ is smooth with integral fibers
 and $F_s$ is torsion free
of rank say $m$ for every $s\in S$.

Over an integral scheme, a torsion free coherent sheaf of rank $m$
is isomorphic to a subsheaf of a free  sheaf of rank $m$.
 Hence there is an exact sequence
 $$
 0\to F|_U\to \o_U^m\to Q\to 0
 $$
 where $\dim Q\leq n-1$.  We also need that $Q$ is $S_{r-1}$,
see for example
\cite[2.60]{kk-singbook}.
In particular, if $r\geq 2$ then $Q$ has pure dimension $n-1$ but if $r=1$
then $Q$ need not be pure dimensional.
 
Applying $j_*$ we get an exact sequence
 $$
 R^{r-1}j_*\o_U^m\to   R^{r-1}j_* Q\to R^rj_*(F|_U)\to R^rj_*\o_U^m\to R^rj_*Q.
  $$
Depending on the values of $(r,n)$, several of the terms vanish.

(\ref{free+comm.main.ptop}.1) If $r>1$ and $n\geq r+2$  then 
$R^{r-1}j_*\o_U=R^rj_*\o_U=0$, thus
$$
R^{r-1}j_* Q\cong R^rj_*(F|_U)
 $$
and we are done by induction.

(\ref{free+comm.main.ptop}.2) If $r>1$ and $n= r+1$  then 
$R^{r-1}j_*\o_U=0$ and $R^rj_*Q=0$, thus we have the exact sequence
 $$
0\to   R^{r-1}j_* Q\to R^rj_*(F|_U)\to R^rj_*\o_U^m\to 0.
  $$
Here $R^{r-1}j_* Q $ and $R^rj_*\o_U^m $ are generically free and 
commute with base change
by induction and  (\ref{free+comms.mooth.lem}).
Thus  $R^rj_*(F|_U) $  is generically free and commutes with base change
by (\ref{free+comm.lem}.1).

(\ref{free+comm.main.ptop}.3) 
If $r=1$ and $n\geq 3$  then $R^1j_*\o_U=0$ and we have the exact sequence
$$
 \o_X^m\to   j_* Q\to R^1j_*(F|_U)\to 0.
  $$
We can use (\ref{free+comm.lem}.2) to show that
$R^1j_*(F|_U) $  is generically free and commutes with base change.

(\ref{free+comm.main.ptop}.4)
 Finally, if $r=1$ and $n=2$ then  we can use the exact sequence
$$
0\to \coker\bigl[ \o_X^m\to   j_* Q\bigr]
\to R^1j_*(F|_U)\to R^1j_*\o_U^m\to 0. \qed
  $$

\section{Restriction of torsion bundles }\label{pf.of.7.sec}

Here we complete the proof of Theorem \ref{gl-gen.thm2.cor}.
We assume that $\depth_xX\geq 3$, hence $\picls(x,X)$ is 0-dimensional
and a 0-dimensional linear algebraic group is finite.
This shows (\ref{gl-gen.thm2.cor}.2--3).
Thus the only new claim is that 
$r^X_D: \picls(x,X){\longrightarrow} \picls(x,D)$ is injective
(\ref{gl-gen.thm2.cor}.1).

We have already noted in Paragraph \ref{unip.last.step} that
every torsion element in the kernel of $r^X_D$
has $p$-power order  where $p$ is the  characteristic. 
The proof relied on the observation that, for $p\nmid m$,
a torsion element of order $m$ in $\picl(x,X) $
corresponds to a degree $m$ \'etale cover of $X\setminus \{x\}$
and then used \cite[XIII.2.1]{sga2}.

Here we develop another approach to deal with
torsion elements in the kernel. The advantage is that 
for this method the
characteristic does not matter. A disadvantage is that
the main step works best for proper morphisms, thus it applies
only when $X$ is essentially of finite type.

\begin{thm} \label{tors.free.ker.thm}
Let  $T$ be the spectrum of a DVR   
with closed point $0\in T$ and generic point $g\in T$.
Let $(x,X)$ be a local scheme and
  $f:(x,X)\to (0,T)$   a flat morphism that is  essentially of finite type.
Let $x\in Z_0\subset X_0$ be a closed subscheme 
such that $\depth_{Z_0}X_0\geq 2$.
Then
$$
\ker\bigl[r^X_{X_0}:\pic(X\setminus Z_0)\to \pic(X_0\setminus Z_0)\bigr]
\qtq{is torsion free.}
$$
\end{thm}

Before we start the proof, we recall some results
on pull-back and push-forward of sheaves.

\begin{say}[Push forward and restriction]
 \label{general.C.crit.ques}
Let $X$ be an excellent scheme, $Z\subset X$  a closed subscheme and
$U:=X\setminus Z$ with natural injection $j:U\into X$.

Let $F$ be a coherent sheaf on $U$ such that 
$\codim_{\bar P}(\bar P\cap Z)\geq 2$ 
for every associated 
prime $P\in U$ of $F$.
Then $j_*F$ is  coherent by \cite[IV.5.11.1]{ega}.
Moreover, it is the unique coherent sheaf $G$ on $X$ such that
$G|_U\cong F$ and $\depth_ZG\geq 2$.

Let  $f:X\to S$  be a flat  morphism  to a regular, 1-dimensional scheme.
For $s\in S$ the restriction of $j$ is denoted by $j_s:U_s\into X_s$. 
There are natural maps
$$
r_s: \bigl(j_*F\bigr)|_{X_s}\to (j_s)_*\bigl(F|_{X_s}\bigr).
\eqno{(\ref{general.C.crit.ques}.1)}
$$
Note that  $j_*F|_{X_s}$  has  $\depth\geq 1$ along $Z\cap X_s$, 
in particular, it has no embedded
points supported in $Z\cap X_s$. 
Thus  $r_s$ is an injection that is an isomorphism over $U_s$.
By \cite[IV.12.1.6]{ega}, $j_*F|_{X_s}$  has  $\depth\geq 2$ 
along $Z\cap X_s$ for general $s\in S$. Thus $r_s$ is
an isomorphisms for general $s\in S$.

We are mainly interested in the case when
  $T=S$ is the spectrum of a DVR   
with closed point $0\in T$ and generic point $g\in T$.
Let  $f:X\to T$ be  a flat morphism and
 $Z_0\subset X_0$  a closed subscheme  such that $\depth_{Z_0}X_0\geq 2$.
This implies that  $\codim_{X_0}Z_0\geq 2$, thus
$\dim X_0\geq 2$  (unless $Z_0$ is empty). 
Set $U:=X\setminus Z_0$ and $U_0:=X_0\setminus Z_0$
with natural open injections  $j:U\into X$ and $j_0:U_0\into X_0$.

Let $L_U$ be a line bundle on $U$ and $L_{U_0}:=L_U|_{U_0}$
its restriction to $U_0$.
We are interested in the  sheaves
$j_*L_U$ and $(j_0)_*L_{U_0}$.

These sheaves are locally free on $U$ (resp.\ $U_0$) and
 have $\depth\geq 2$ along $Z_0$. 
If $L_{U_0}\cong \o_{U_0}$ then  $(j_0)_*L_{U_0}\cong \o_{X_0}$;
here we use the assumption 
 that  $\depth_{Z_0}X_0\geq 2$.

As in (\ref{general.C.crit.ques}.1) there is  a natural map 
$$
r_0:j_*L_U|_{X_0}\to (j_0)_*L_{U_0}.
\eqno{(\ref{general.C.crit.ques}.2)}
$$
Note that $r_0$ an injection that is an isomorphism over $U_0$.
We aim to understand when it is an isomorphism along $Z_0$
and then relate this to (\ref{tors.free.ker.thm}) through a
series of local freeness criteria for $j_*L_U$.
\end{say}

\begin{lem} \label{torsion.ker.in.r.lem.0}
Notation and assumptions as in (\ref{general.C.crit.ques}).
Then  $j_*L_U$ is locally free 
iff 
\begin{enumerate}
\item $(j_0)_*L_{U_0}$ is locally free  and 
\item $r_0:j_*L_U|_{X_0}\to (j_0)_*L_{U_0}$ is an isomorphism.
\end{enumerate}
\end{lem}

Proof. If  $j_*L_U$ is locally free  then $j_*L_U|_{X_0}$ is locally free
hence it has  $\depth\geq 2$ along $Z_0$. Thus
 $(j_0)_*L_{U_0}=j_*L_U|_{X_0}$ is locally free.

Conversely, if $j_*L_U|_{X_0}$ is locally free then
$j_*L_U$ is locally free near $X_0$ by the Nakayama lemma. \qed

\begin{say}[Numerical inequalities]\label{torsion.ker.in.r.lem.1.num}
Notation and assumptions as in (\ref{general.C.crit.ques}).
Assume in addition that $f$ is proper and has relative dimension $n$.
By semicontinuity
$$
H^0\bigl(X_0, (j_0)_*L_{U_0}\bigr)\geq  H^0\bigl(X_0, j_*L_U|_{X_0}\bigr)\geq  
H^0\bigl(X_g, L_g\bigr)
\eqno{(\ref{torsion.ker.in.r.lem.1.num}.1)}
$$
where $L_g:=L_U|_{X_g}$. In particular, if $L_g$ is ample on $X_g$
and $(j_0)_*L_{U_0}$ is locally free and ample on $X_0$
then, applying (\ref{torsion.ker.in.r.lem.1.num}.1)
to powers of $L_U$ we conclude that
$$
\bigl(c_1((j_0)_*L_{U_0})^n\bigr)\geq \bigl(c_1(L_g)^n\bigr).
\eqno{(\ref{torsion.ker.in.r.lem.1.num}.2)}
$$
One can be more precise if  $\dim Z_0=0$. 
Then  $\coker r_0 $ is artinian, thus
$$
\chi\bigl(X_0, (j_0)_*L_{U_0}\bigr)=
\chi\bigl(X_0, j_*L_U|_{X_0}\bigr)+\len(\coker r_0).
\eqno{(\ref{torsion.ker.in.r.lem.1.num}.3)}
$$
Since $j_*L_U $ is flat over $T$,
$$
\chi\bigl(X_0, j_*L_U|_{X_0}\bigr)=
\chi\bigl(X_g, j_*L_U|_{X_g}\bigr)=
\chi\bigl(X_g, L_g\bigr).
$$
Combining these we see that 
$$
\chi\bigl(X_0, (j_0)_*L_{U_0}\bigr)=
\chi\bigl(X_g, L_g\bigr)+\len(\coker r_0).
\eqno{(\ref{torsion.ker.in.r.lem.1.num}.4)}
$$
\end{say}

\begin{cor} \label{torsion.ker.in.r.lem.1}
Notation and assumptions as in (\ref{general.C.crit.ques}).
Assume in addition that $f$ is proper and $\dim Z_0=0$. 
Then  $j_*L_U$ is locally free 
iff 
\begin{enumerate}
\item $(j_0)_*L_{U_0}$ is locally free  and 
\item $\chi\bigl(X_0, (j_0)_*L_{U_0}\bigr)=\chi\bigl(X_g, L_g\bigr)$.
\end{enumerate}
\end{cor}

Proof.  By (\ref{torsion.ker.in.r.lem.1.num}.4)
$\chi\bigl(X_0, (j_0)_*L_{U_0}\bigr)=
\chi\bigl(X_g, L_g\bigr)+\len(\coker r_0)$.
Thus $r_0$ is an isomorphism if
$\chi\bigl(X_0, (j_0)_*L_{U_0}\bigr)=\chi\bigl(X_g, L_g\bigr)$
hence (\ref{torsion.ker.in.r.lem.0}) applies.\qed

\medskip

The following is a global version of Theorem \ref{tors.free.ker.thm}
and a key step of its proof.

\begin{lem} \label{torsion.ker.in.r.lem.2}
Notation and assumptions as in (\ref{general.C.crit.ques}).
Assume in addition that $f$ is proper and $\dim Z_0=0$. 
Then  $j_*L_U$ is locally free  iff 
\begin{enumerate}
\item $(j_0)_*L_{U_0}$ is locally free  and 
\item $j_*\bigl(L_U^r\bigr)$ is locally free  for some $r>0$.
\end{enumerate}
\end{lem}

Proof.  Note that if  $(j_0)_*L_{U_0}$ is locally free then so is
$$
\bigl((j_0)_*L_{U_0}\bigr)^m\cong (j_0)_*\bigl(L_{U_0}^m\bigr)
\qtq{for every $m\in \z$.}
$$
We use (\ref{torsion.ker.in.r.lem.1}) for $L_U^m$ for every $m$.
We obtain that the following are equivalent:
\begin{enumerate}\setcounter{enumi}{2}
\item $j_*L_U$ is locally free, 
\item $j_*\bigl(L_U^m\bigr)$ is locally free for every $m\in \z$ and
\item
$\chi\bigl(X_0, \bigl((j_0)_*L_{U_0}\bigr)^m\bigr)=
\chi\bigl(X_g, L_g^m\bigr)$ for every $m\in \z$.
\end{enumerate}
Both Euler characteristics are polynomials in $m$, hence they agree
iff they agree for infinitely many values of $m$.
Since   $j_*\bigl(L_U^r\bigr)$ is locally free we know that
$$
\chi\bigl(X_0, \bigl((j_0)_*L_{U_0}\bigr)^{rm}\bigr)=
\chi\bigl(X_g, L_g^{rm}\bigr)
\qtq{for every $m\in \z$.}
$$
Thus $\chi\bigl(X_0, (j_0)_*L_{U_0}\bigr)=
\chi\bigl(X_g, L_g\bigr)$
and so $j_*L_U$ is locally free by (\ref{torsion.ker.in.r.lem.1}).\qed
\medskip

Next we remove the properness assumption
from (\ref{torsion.ker.in.r.lem.2}).

\begin{lem} \label{torsion.ker.in.r.lem.2.aff}
Notation and assumptions as in (\ref{general.C.crit.ques}).
Assume in addition that $f$ is of finite type and $\dim Z_0=0$. 
Then  $j_*L_U$ is locally free  iff 
\begin{enumerate}
\item $(j_0)_*L_{U_0}$ is locally free  and 
\item $j_*\bigl(L_U^r\bigr)$ is locally free  for some $r>0$.
\end{enumerate}
\end{lem}

Proof. The assertions are local on $X$, hence we may
assume that $X$ is affine.  Since  $f$ is of finite type,
there is a compactification  $\bar X\supset X$ such that
$f$ extends to a  flat morphism $\bar f:\bar X\to T$.
We intend to use Lemma \ref{torsion.ker.in.r.lem.2},
the only problem is that we do not know 
how to extend $L_U$ to $\bar X$ and
what happens along
 $\bar X\setminus X$. 

Thus we need to improve  the compactification $\bar X$.
In order to do this,
choose an effective  Cartier divisor
$D_U$ on $U$ such that $L_U\cong \o_U(-D_U)$.
Let $\bar D\subset \bar X$ denote its closure
and  $I_{\bar D}:= \o_{\bar X}(-\bar D)$ 
its ideal sheaf. 
Let  $\tilde X$ be the scheme obtained by gluing
$\tilde X_1:=X$ and the blow-up  
$$
\tilde X_2:=B_{I_{\bar D}}(\bar X\setminus Z_0)
$$
along  $\tilde X_1\supset U\cong B_{I_{\bar D}}U\subset \tilde X_2$. 
By construction $f$ extends to  a  flat morphism $\tilde f:\tilde X\to T$.
Note that $\tilde f$ is proper  (but not necessarily projective).
Also, even if $f$ has $S_2$ fibers, the
central fiber of  $\tilde f$ need not be $S_2$. 
Let $\tilde X_0$ denote the fiber of $\tilde f$ over $0\in T$. 
Set  $\tilde U:=\tilde X\setminus Z_0$ and
 $\tilde U_0:=\tilde X_0\setminus Z_0$
with natural injections
$\tilde j:\tilde U\into \tilde X$ and
$\tilde j_0:\tilde U_0\into \tilde X_0$.

Let $I_{\tilde D}\subset \o_{\tilde X}$ denote the inverse image ideal sheaf of
$I_{\bar D}$.
 Then $I_{\tilde D}$ is locally free on $\tilde U:=\tilde X\setminus Z_0$;
denote its restriction by $L_{\tilde U}$.
By construction  $L_{\tilde U}|_U\cong L_U$ and hence
 $\tilde j_*L_{\tilde U}^r$ is locally free.
Similarly, $L_{\tilde U_0}$ agrees with $L_{U_0}$
over $U_0\subset X_0\subset \tilde X_0$, thus
$(\tilde j_0)_*L_{\tilde U_0}$ is locally free.

Thus Lemma \ref{torsion.ker.in.r.lem.2}
applies to $\tilde f:\tilde X\to T$ and  $L_{\tilde U}$
to conclude that $\tilde j_*L_{\tilde U}$ is locally free.
Therefore $j_*L_U$ is also locally free.
\qed

 \begin{say}[Proof of Theorem \ref{tors.free.ker.thm}]
 \label{tors.free.ker.thm.pf}
We use the notation of (\ref{general.C.crit.ques}). 
Let $L_U$ be a line bundle on $U$ such that  $L_{U_0}:=L_U|_{U_0}\cong \o_{U_0}$
and $L_U^r\cong \o_U$ for some $r>0$. 

If $r_0: j_*L_U|_{X_0}\to (j_0)_*L_{U_0}$ is an isomorphism then
the constant 1 section of  $L_{U_0}\cong \o_{U_0}$ lifts to
a nowhere zero section of $L_U$ and so $L_U\cong \o_{U}$.
Thus it is enough to prove that  
$r_0: j_*L_U|_{X_0}\to (j_0)_*L_{U_0}$ is an isomorphism.
The  latter  is a local question on $X$. 

By localizing at a generic point of the set where $r_0$ is not known to be
an isomorphism we are reduced to the case when $Z_0$ is a single closed point
of $X_0$. As in Definition
\ref{loc.pic.alg.defn},  we are further reduced to the case when
$f:X\to T$ is of finite type and the latter was treated in
Lemma \ref{torsion.ker.in.r.lem.2.aff}.\qed
\end{say}

\section{Numerical criteria for relative Cartier  divisors}
\label{num.crit.sec}

\begin{defn}\label{gen.flat.divs.defn}
 Let  $T$ be a regular 1-dimensional scheme and
 $f:X\to T$  a flat morphism  with  $S_2$ fibers.
A {\it generically flat family of divisors}  $D$ on $X$ 
is given by
\begin{enumerate}
\item  an open set $U\subset X$ such that
  $\codim_{X_t}\bigl(X_t\setminus U\bigr)\geq 2$ for every $t\in T$ and
\item a relative Cartier divisor $D_U$ on $U$; that is, a Cartier divisor
whose restriction to every fiber is a  Cartier divisor.
\end{enumerate}
For each $t\in T$ set $U_t:=X_t\cap U$. Then 
 $D_U|_{U_t}$ extends uniquely to a  divisor 
 on $X_t$; we denote it by $D_t$ and call it the {\it restriction}
 of $D$ to $X_t$.

If $|H|$ is a base point free linear system on $X$ then
the restriction of $D$ to a general $H\in |H|$
is again generically flat family of divisors.

We say that $D$ is a {\it fiber-wise Cartier} family of divisors
if each $D_t$ is a Cartier divisor. 
As we noted in (\ref{general.C.crit.ques}), 
$$
r_t: \o_X(D)|_{X_t}=\bigl(j_*\o_U(D_U)\bigr)|_{X_t}\to 
(j_t)_*\bigl(\o_{U_t}(D_{U_t})\bigr)=\o_{X_t}(D_t)
$$
is an isomorphism for general $t\in T$; thus $D$ is Cartier 
except possibly along finitely many closed fibers.

Our main interest is to find conditions 
that guarantee that a fiber-wise Cartier family of divisors is  everywhere
Cartier.

We say that $D$ is a {\it fiber-wise ample} family of divisors
if each $D_t$ is an ample Cartier divisor.
\end{defn}

Combining   Theorems  \ref{gl-gen.thm2.cor}
and  \ref{tors.free.ker.thm} we have the following.

 \begin{thm}\label{num.C.crit.ques.codim3}
Let  $T$ be an irreducible,  regular, 1-dimensional scheme and
 $f:X\to T$  a flat  morphism that is essentially of finite type and has
  $S_2$ fibers.

Let $D$ be a generically flat family of fiber-wise Cartier divisors on $X$.
Assume that  there is a closed subscheme  $W\subset X$ such that
 $\codim_{X_t}\bigl(X_t\cap W\bigr)\geq 3$ for every $t\in T$ and
$D$ is   Cartier  on $U:=X\setminus W$.

Then  $D$ is a  Cartier divisor  on $X$. \qed
\end{thm}

It remains to understand what happens when the putative non-Cartier locus
has codimension 2.
For  projective morphisms we have the following numerical criteria.

 \begin{thm}\label{num.C.crit.ques}
Let  $T$ be an irreducible,  regular, 1-dimensional scheme and
 $f:X\to T$  a flat, projective morphism 
of relative dimension $n$
with  $S_2$ fibers.

Let $D$ be a generically flat family of fiber-wise Cartier divisors on $X$
and $H$  an $f$-ample Cartier divisor on $X$. The following are equivalent.
\begin{enumerate}
\item $D$ is a Cartier divisor on $X$.
\item $\bigl(D_t^2\cdot H_t^{n-2}\bigr)$ is independent of $t\in T$.
\end{enumerate}
\noindent If $D$ is fiber-wise ample, then these are further  equivalent to
\begin{enumerate}\setcounter{enumi}{2}
\item 
$\bigl(D_t^n\bigr)$ is independent of
$t\in T$.
\end{enumerate}
\end{thm}

Proof. If $D$ is Cartier then all  the intersection numbers
$\bigl(D_t^i\cdot H_t^{n-i}\bigr)$ are independent of
$t\in T$. Thus (1) $\Rightarrow$ (2) and (3). 

To see the converse we may assume that
$T$ is local with closed point $0\in T$ and generic point $g\in T$.
Let $Z_0\subset X$ be the smallest closed subset 
such that $D$ is Cartier on $X\setminus Z_0$. 
Note that $Z_0\subset X_0$ since $D$ is Cartier on $X_g$.
We can choose
$U:= X\setminus Z_0$ as the open set in (\ref{gen.flat.divs.defn}).

Assume first that $n=2$, thus (2) and (3) coincide.
For each $t\in T$, the  Euler characteristic
is a quadratic polynomial
$$
\chi\bigl(X_t, \o_{X_t}(mD_t)\bigr)=a_tm^2+b_tm+c_t,
$$
and we know from Riemann--Roch that 
$a_t=\tfrac12\bigl(D_t^2\bigr)$ and $c_t=\chi\bigl(X_t, \o_{X_t}\bigr)$.
Furthermore,  (\ref{torsion.ker.in.r.lem.1.num}.4) implies that
$$
a_0m^2+b_0m+c_0\geq a_gm^2+b_gm+c_g\qtq{for every $m\in \z$.}
\eqno{(\ref{num.C.crit.ques}.4)}
$$ 
For $m\gg 1$ the quadratic terms dominate, which gives that
$$
\bigl(D_0^2\bigr)=2a_0\geq 2a_g=\bigl(D_g^2\bigr).
\eqno{(\ref{num.C.crit.ques}.5)}
$$
Assume now that $(D_0^2)=(D_g^2 ) $. Then $a_0=a_g$ thus 
(\ref{num.C.crit.ques}.5) implies that
$$
b_0m+c_0\geq b_gm+c_g\qtq{for every $m\in \z$.}
\eqno{(\ref{num.C.crit.ques}.6)}
$$
For $m\gg 1$ this implies that $b_0\geq b_g$ and for
$m\ll -1$ that $-b_0\geq -b_g$.
Thus $b_0=b_g$ and  
$c_0=\chi\bigl(X_0, \o_{X_0}\bigr)=\chi\bigl(X_g, \o_{X_g}\bigr)=c_g$  
also holds since $f$ is flat. 
Therefore we have equality in (\ref{num.C.crit.ques}.4). 
Using (\ref{torsion.ker.in.r.lem.1}) we see that
$$
r_0: j_*\o_U(D_U)|_{X_0}\to (j_0)_*\o_{U_0}(D_0)
$$
 is an isomorphism and  hence $\o_X(D)=j_*\o_U(D_U)$ is locally free.

In order to prove (2) $\Rightarrow$ (1)
for $n\geq 3$ we use induction on $n$.
A suitable multiple $|mH|$ provides an embedding  $X\subset \p^N_T$; let
$X'\subset X$ be a general hyperplane section. 
Then  $f':=f|_{X'}:X'\to T$  is a flat, projective morphism 
of relative dimension $n-1$
 and  $Z':=Z\cap X'$ is a closed subscheme such that
 $\codim_{X'_t}\bigl(X'_t\cap Z'\bigr)\geq 2$ for every $t\in T$.

Furthermore, $D':=D|_{X'}$ is a divisor, 
 Cartier over  $U':=X'\setminus Z'$ whose
 support does not contain any irreducible component of a fiber of $f'$
and   $D'|_{X'_t}$ is   Cartier  for every $t$. Finally
$\bigl((D'_t)^2\cdot H_t^{n-3}\bigr)=m\bigl(D_t^2\cdot H_t^{n-2}\bigr)$
is independent of $t\in T$.

Since $f$ has  $S_2$ fibers, a hyperplane section of it usually has only
$S_1$ fibers. However, by
\cite[IV.12.1.6]{ega}, a general hyperplane section 
again has $S_2$ fibers.
Thus, by induction, $D'=D|_{X'}$ is a  Cartier divisor.
As we noted in (\ref{gen.flat.divs.defn}), since $X'$ is general,
this implies that  $D$ is Cartier along $X'$. Hence
there is a closed subscheme $W\subset X$ such that
$W\cap X_t$ is $0$-dimensional for every $t\in T$
and $D$ is  Cartier on $X\setminus W$. 
Thus $D$ is Cartier by  Theorem \ref{num.C.crit.ques.codim3}.

Assume finally that $D$ is fiber-wise ample. Choose $m>0$
such that $mD_0-H_0$ and $mD_t-H_t$ are ample. Thus
$B:=mD-H$ is also fiber-wise ample and 
$$
m^n\bigl(D_t^n\bigr)=\tsum_{i=0}^n \bigl(B_t^i\cdot H_t^{n-i}\bigr).
$$
We can apply (\ref{torsion.ker.in.r.lem.1.num}.2)
to the restriction of $B$ to an intersection of $n-i$ general members of
$H$ to obtain that 
$\bigl(B_0^i\cdot H_0^{n-i}\bigr)\geq \bigl(B_g^i\cdot H_g^{n-i}\bigr)$
for every $i$. However, $\bigl(D_0^n\bigr)= \bigl(D_g^n\bigr)$
by assumption (3), thus
$\bigl(B_0^i\cdot H_0^{n-i}\bigr)= \bigl(B_g^i\cdot H_g^{n-i}\bigr)$ for every $i$.
We can use this for $i=2$ and the already established (2) $\Rightarrow$ (1)
to conclude that $B$ is Cartier. Thus $mD=B+H$ is also Cartier.
By the same argument $(m+1)D$ is also Cartier;
 this  implies that $D$ is Cartier.
\qed

\section{Open problems}

\subsection*{Finite generation of $\nsl(x,X)$}{\ }

The main unsolved problem is the finite generation of $\nsl(x,X)$.
This is known if $X$ is normal or if the characteristic is 0 and $\dim X\geq 3$.
In positive characteristic, our methods do not seem to distinguish
a unipotent subgroup of $\piclo$ from a discrete $p$-group.

For surfaces over $\c$, even stronger results should hold; see
(\ref{pic.H2.lem.conj}). 

\subsection*{Local Picard group of excellent schemes}{\ }

While our theorems settle only the geometric cases, I see no reason
why they should not hold in general.

\begin{conj}\label{pic-norm-loc.exc.conj}
Let  $X$ be an excellent scheme  that is $S_2$
 and has pure dimension $\geq 3$.
Let $x\in X$ be a closed point.
Then 
$ \ker\bigl[\pi^*: \picls(x,X)\to
 \picls(\bar x,\bar X) \bigr]$
is a linear algebraic group.
\end{conj}

\begin{conj} \label{gl-gen.thm2.conj}
Let  $X$ be an excellent scheme  that is $S_2$
 and has pure dimension $\geq 4$. Let $x\in X$ be a closed point
 and  $x\in D\subset X$  an effective Cartier divisor. 
Then  
$ \ker\bigl[r^X_D:\picls(x,X)\to
 \picls(x, D) \bigr]$
is a unipotent algebraic group.
\end{conj}

Note that not even the existence of $\picls(x,X) $ 
is known in the above generality, but algebraic equivalence
can be defined as in (\ref{loc.pic.func.defn}.a--c).

\subsection*{Local Picard group of surfaces}{\ }

The construction of \cite{MR492263} does not apply
to surfaces; in fact, the functorial approach always  gives
the ``wrong'' answer. For normal surface singularities
 over $\c$, \cite{mumf-top} constructed
a finite dimensional  local Picard group  $\picls(s,S)$.
The examples (\ref{pic-norm-loc.exmp.1}.3--4) show that in many non-normal cases
there is no finite dimensional  local Picard group. 

\begin{conj} Let $(s, S)$ be an $S_2$ surface. Then
there is a  finite dimensional  local Picard group
$\picls(s, S)$  iff $S$ is seminormal.
\end{conj}

\begin{conj} \label{pic.H2.lem.conj} Let $(x, X$) be a local 
$\c$-scheme of finite type
that is $S_2$ and has pure dimension $2$.
Then taking the first Chern class gives an exact sequence
$$
0\to \piclo(x,X)\to \picls(x,X)\stackrel{c_1}{\to}
H^2\bigl(\link(x,X), \z\bigr).
$$
\end{conj}

The normal case is discussed in \cite{mumf-top}.
It would be nice to have something similar in
positive characteristic, but I do not know what
should replace  $H^2\bigl(\link(x,X), \z\bigr) $.
See also Example \ref{pinch.exmp}.

In connection with (\ref{pic.not.glob.exmp}), one can ask
the following.

\begin{ques} Let $(0, S^{\rm an})$ be a normal, analytic
 surface singularity over $\c$
such that $\piclo(0, S)$ is proper. Is there an algebraic model
$(0, S)$ such that $\picl(0, S)= \picl(0, S^{\rm an})$?

\end{ques}

\subsection*{Second cohomology of links}{\ }

One can refine the topological approach 
of Section \ref{anal.sect}
using the
mixed Hodge structures on the 
cohomology groups of the links $H^2\bigl(\link(x,X), \c\bigr)$
 and $H^2\bigl(\link(\bar x,\bar X), \c\bigr)$
\cite[Sec.6.2]{PetersSteenbrinkBook}.
Since the Chern class of a line bundle
has pure Hodge type $(1,1)$, the following 
would imply most of
Conjecture \ref{pic-norm-loc.exc.conj} for analytic spaces.

\begin{conj}\label{pic-norm-loc.anal-wt2.conj}
Let  $X$ be a complex analytic space that is  $S_2$ and has
pure dimension $\geq 3$.
Let $x\in X$ be a point.
Then  the pull-back map
$$
\pi^* : H^2\bigl(\link(x,X), \c\bigr)\to H^2\bigl(\link(\bar x,\bar X), \c\bigr)
$$
is injective on the weight 2 graded piece of the mixed Hodge structure.
\end{conj}

\subsection*{Obstruction theory for $\picls$}{\ }

The usual obstruction space for  $\picls(x,X)$ is
$H^3_x(X, \o_X)$. This is infinite dimensional if $\dim X=3$ or if
$\dim X\geq 4$ and  $X\setminus \{x\}$ is not $S_3$.

Already  \cite{Artin74c} observed that an
obstruction  theory is an external construct imposed on a functor
and there could be different ``natural'' 
obstruction  theories that work for the same
functor.
 Theorem \ref{free+comm.main.ptop} suggests
that one should be able to develop a finite dimensional obstruction theory
for $\picls$.  I do not know how to formulate such a theory.


 \begin{ack}
I thank B.~Bhatt, J.~de~Jong and K.~Smith
 for  answering  my questions and
suggesting many improvements. 
Partial financial support   was provided  by  the NSF under grant number 
DMS-0968337.
\end{ack}

\def\cprime{$'$} \def\cprime{$'$} \def\cprime{$'$} \def\cprime{$'$}
  \def\cprime{$'$} \def\cprime{$'$} \def\dbar{\leavevmode\hbox to
  0pt{\hskip.2ex \accent"16\hss}d} \def\cprime{$'$} \def\cprime{$'$}
  \def\polhk#1{\setbox0=\hbox{#1}{\ooalign{\hidewidth
  \lower1.5ex\hbox{`}\hidewidth\crcr\unhbox0}}} \def\cprime{$'$}
  \def\cprime{$'$} \def\cprime{$'$} \def\cprime{$'$}
  \def\polhk#1{\setbox0=\hbox{#1}{\ooalign{\hidewidth
  \lower1.5ex\hbox{`}\hidewidth\crcr\unhbox0}}} \def\cdprime{$''$}
  \def\cprime{$'$} \def\cprime{$'$} \def\cprime{$'$} \def\cprime{$'$}
\providecommand{\bysame}{\leavevmode\hbox to3em{\hrulefill}\thinspace}
\providecommand{\MR}{\relax\ifhmode\unskip\space\fi MR }
\providecommand{\MRhref}[2]{%
  \href{http://www.ams.org/mathscinet-getitem?mr=#1}{#2}
}
\providecommand{\href}[2]{#2}

\bigskip

\noindent Princeton University, Princeton NJ 08544-1000

{\begin{verbatim}kollar@math.princeton.edu\end{verbatim}}

\end{document}